\documentclass[12pt,reqno]{amsart}

\usepackage{amssymb}

\usepackage{mathrsfs}

\DeclareMathAlphabet{\mathpzc}{OT1}{pzc}{m}{it}

\usepackage[all]{xy}

\entrymodifiers={+!!<0pt,\fontdimen22\textfont2>}

\def\cA{\mathscr{A}}

\def\cC{\mathscr{C}}

\def\cH{\mathscr{H}}

\def\BZ{\mathbb{Z}}

\def\fA{\mathfrak{A}}

\def\fa{\mathfrak{a}}

\def\fx{\mathfrak{x}}

\def\sC{\mathsf{C}}
\def\sD{\mathsf{D}}

\def\add{\operatorname{add}}

\def\adots{\mathinner{\mkern1mu\raise1.0pt\vbox{\kern7.0pt\hbox{.}}\mkern2mu\raise4.0pt\hbox{.}\mkern2mu\raise7.0pt\hbox{.}\mkern1mu}}

\def\C{\operatorname{C}}

\def\D{\sD}

\def\Df{\D^{\operatorname{f}}}

\def\H{\operatorname{H}}

\def\Hom{\operatorname{Hom}}

\def\mod{\mathsf{mod}}

\def\prod{\operatorname{prod}}

\def\rep{\operatorname{rep}}

\newtheorem{Lemma}{Lemma}[section]
\newtheorem{Theorem}[Lemma]{Theorem}
\newtheorem{Proposition}[Lemma]{Proposition}
\newtheorem{Corollary}[Lemma]{Corollary}

\theoremstyle{definition}
\newtheorem{Definition}[Lemma]{Definition}
\newtheorem{Setup}[Lemma]{Setup}

\newtheorem{Remark}[Lemma]{Remark}

\newtheorem{Example}[Lemma]{Example}

\newtheorem{Notation}[Lemma]{Notation}

\def\arcs{\fA}
\def\cluster{A}
\def\subcat{\cA}

\begin{document}

\setlength{\parindent}{0pt}
\setlength{\parskip}{7pt}
%The default \baselineskip is close to 4.8mm
%\setlength{\baselineskip}{5.3mm}

\title[Cluster category of infinite Dynkin type]
{On a triangulated category which behaves like a cluster category of
infinite Dynkin type, and the relation to triangulations of the
infinity-gon}

\author{Thorsten Holm}
\address{Institut f\"{u}r Algebra, Zahlentheorie und Diskrete
  Mathematik, Fa\-kul\-t\"{a}t f\"{u}r Ma\-the\-ma\-tik und Physik, Leibniz
  Universit\"{a}t Hannover, Welfengarten 1, 30167 Hannover, Germany}
%{Institut f\"{u}r Algebra und Geometrie, Fakult\"{a}t f\"{u}r
%Mathematik, Otto-von-Guericke-Universit\"{a}t Magdeburg, Postfach
%4120, 39016 Magdeburg, Germany\newline
%and\newline
%Department of Pure Mathematics, University of Leeds,
%Leeds LS2 9JT, United Kingdom}
\email{holm@math.uni-hannover.de}
\urladdr{http://www.iazd.uni-hannover.de/\~{ }tholm}

\author{Peter J\o rgensen}
\address{School of Mathematics and Statistics,
Newcastle University, Newcastle upon Tyne NE1 7RU,
United Kingdom}
\email{peter.jorgensen@ncl.ac.uk}
\urladdr{http://www.staff.ncl.ac.uk/peter.jorgensen}

%\author{Next author goes here}
%\address{Next author's postal address goes here}
%\email{Next author's mail address goes here}

%\thanks{Date: \today. A thank you would go here}


%16E10: Homological dimension
%16E45: Differential graded algebras and applications
%16G10: Representations of Artinian rings 
%16G20: Representations of quivers and partially ordered sets
%16G60: Representation type (finite, tame, wild, etc.) 
%16G70: Auslander-Reiten sequences (almost split sequences) and
%       Auslander-Reiten quivers
%18E30: Derived categories, triangulated categories
%18E35: Localization of categories 
%18G99: Homological algebra: None of the above, but in this section 
%55P62: Rational homotopy theory

% \begin{abstract} 

% We investigate a triangulated category which behaves like a cluster
% category of type $A_{\infty}$.
  
% \end{abstract}

\maketitle

\setcounter{section}{-1}
\section{Introduction}
\label{sec:introduction}

This paper investigates a certain $2$-Calabi-Yau triangulated category
$\sD$ whose Aus\-lan\-der-Rei\-ten quiver is $\BZ A_{\infty}$.  We
show that the cluster tilting subcategories of $\sD$ form a so-called
cluster structure, and we classify these subcategories in terms of
what one may call `triangulations of the $\infty$-gon'.

This is reminiscent of the cluster category $\sC$ of type $A_n$ which
is a $2$-Calabi-Yau triangulated category whose Auslander-Reiten
quiver is a quotient of $\BZ A_n$; see \cite{BMRRT} and \cite{CCS}.
The cluster tilting subcategories of $\sC$ form a cluster structure
and they are classified in terms of triangulations of the $(n+3)$-gon
by
\cite{CCS}.

The category $\sD$ behaves like a `cluster category of type
$A_{\infty}$'.

Let us give some more details.  The category $\sD$ is the compact
derived category $\sD^c(A)$ of the Differential Graded cochain algebra
$A = \C^*(S^2;k)$ where $S^2$ is the $2$-sphere and $k$ is a field.
In fact, $S^2$ is formal over any field $k$, so one can equally well
use as $A$ the quasi-isomorphic DG algebra obtained by placing copies
of $k$ in cohomological degrees $0$ and $2$.

%The AR quiver of $\sD$ is $\BZ A_{\infty}$; see \cite[sec.\ 8]{J} and
%\cite[thm.\ 1.8]{K}.

The category $\sD$ was studied in \cite{J}, and it follows from
\cite[prop.\ 4.4 and cor.\ 5.2]{J} that it is a $2$-Calabi-Yau
category.  This makes it natural to consider maximal $1$-orthogonal
subcategories $\subcat$ of $\sD$.  These are the subcategories which
satisfy $\subcat = (\Sigma^{-1}\subcat)^{\perp} = {}^{\perp}(\Sigma
\subcat)$ and they were introduced by Iyama; see \cite{BMRRT},
\cite{Iyama2}, \cite{Iyama1}, \cite{IyamaYoshino}, and
\cite{KellerReiten2}.  Our first main result is this.

{\bf Theorem A. }
{\em
There is a bijection between maximal $1$-orthogonal subcategories of 
$\sD$ and triangulations of the $\infty$-gon.
}

By a triangulation of the $\infty$-gon, we mean a maximal set of
non-in\-ter\-sec\-ting arcs connecting non-neighbouring integers: We
adopt the philosophy that the integers can be viewed as the vertices
of the $\infty$-gon, and that the arcs can be viewed as diagonals.
There are two obvious ways to achieve such maximal sets; they are
shown in the following two sketches where the arcs must be continued
ad infinitum according to the indicated pattern.  First a `leapfrog'
configuration which is locally finite in the sense that only finitely
many arcs end in each integer.
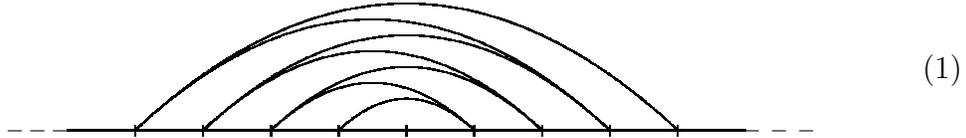
\begin{equation}
\label{equ:leapfrog}
\vcenter{
  \xymatrix @-4.75pc @! {
       \rule{0ex}{10ex} \ar@{--}[r]
     & *{}\ar@{-}[r]
     & *{\rule{0.1ex}{0.8ex}} \ar@{-}[r] \ar@/^3.5pc/@{-}[rrrrrrr]\ar@/^4.0pc/@{-}[rrrrrrrr]
     & *{\rule{0.1ex}{0.8ex}} \ar@{-}[r] \ar@/^2.5pc/@{-}[rrrrr]\ar@/^3.0pc/@{-}[rrrrrr]
     & *{\rule{0.1ex}{0.8ex}} \ar@{-}[r] \ar@/^1.5pc/@{-}[rrr]\ar@/^2.0pc/@{-}[rrrr]
     & *{\rule{0.1ex}{0.8ex}} \ar@{-}[r] \ar@/^1.0pc/@{-}[rr]
     & *{\rule{0.1ex}{0.8ex}} \ar@{-}[r] 
     & *{\rule{0.1ex}{0.8ex}} \ar@{-}[r]
     & *{\rule{0.1ex}{0.8ex}} \ar@{-}[r]
     & *{\rule{0.1ex}{0.8ex}} \ar@{-}[r]
     & *{\rule{0.1ex}{0.8ex}} \ar@{-}[r]
     & *{}\ar@{--}[r]
     & *{}
                    }
}
\end{equation}
Then a `fountain' where infinitely many arcs going to either side end
in a single integer.
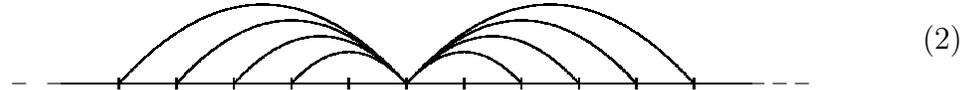
\begin{equation}
\label{equ:fountain}
\vcenter{
  \xymatrix @-4.0pc @! {
       \rule{0ex}{7.5ex} \ar@{--}[r]
     & *{}\ar@{-}[r]
     & *{\rule{0.1ex}{0.8ex}} \ar@{-}[r] \ar@/^2.5pc/@{-}[rrrrr]
     & *{\rule{0.1ex}{0.8ex}} \ar@{-}[r] \ar@/^2.0pc/@{-}[rrrr]
     & *{\rule{0.1ex}{0.8ex}} \ar@{-}[r] \ar@/^1.5pc/@{-}[rrr]
     & *{\rule{0.1ex}{0.8ex}} \ar@{-}[r] \ar@/^1.0pc/@{-}[rr]
     & *{\rule{0.1ex}{0.8ex}} \ar@{-}[r] 
     & *{\rule{0.1ex}{0.8ex}} \ar@{-}[r] 
     & *{\rule{0.1ex}{0.8ex}} \ar@{-}[r]
     & *{\rule{0.1ex}{0.8ex}} \ar@{-}[r] \ar@/^-1.0pc/@{-}[ll]
     & *{\rule{0.1ex}{0.8ex}} \ar@{-}[r] \ar@/^-1.5pc/@{-}[lll]
     & *{\rule{0.1ex}{0.8ex}} \ar@{-}[r] \ar@/^-2.0pc/@{-}[llll]
     & *{\rule{0.1ex}{0.8ex}} \ar@{-}[r] \ar@/^-2.5pc/@{-}[lllll]
     & *{}\ar@{--}[r]
     & *{}
                    }
}
\end{equation}

Maximal $1$-orthogonal subcategories are particularly important if
they are functorially finite.  Then they are called cluster tilting
subcategories and the corresponding quotient categories are abelian by
\cite[sec.\ 2]{KellerReiten2} and \cite[thm.\ 3.3]{KoenigZhu}.  Our
second main result is the following.

{\bf Theorem B. }
{\em
A maximal $1$-orthogonal subcategory of $\sD$ is functorially finite
if and only if the corresponding triangulation of the $\infty$-gon is
locally finite or has a fountain.
}

The point is that there are triangulations of the $\infty$-gon like
the following, which have a `right-fountain' and a `left-fountain' but
do not satisfy the conditions of Theorem B.
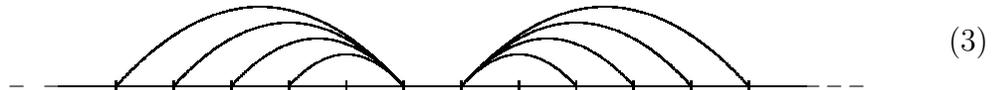
\begin{equation}
\label{equ:right_and_left_fountain}
\vcenter{
  \xymatrix @-4.0pc @! {
       \rule{0ex}{7.5ex} \ar@{--}[r]
     & *{}\ar@{-}[r]
     & *{\rule{0.1ex}{0.8ex}} \ar@{-}[r] \ar@/^2.5pc/@{-}[rrrrr]
     & *{\rule{0.1ex}{0.8ex}} \ar@{-}[r] \ar@/^2.0pc/@{-}[rrrr]
     & *{\rule{0.1ex}{0.8ex}} \ar@{-}[r] \ar@/^1.5pc/@{-}[rrr]
     & *{\rule{0.1ex}{0.8ex}} \ar@{-}[r] \ar@/^1.0pc/@{-}[rr]
     & *{\rule{0.1ex}{0.8ex}} \ar@{-}[r] 
     & *{\rule{0.1ex}{0.8ex}} \ar@{-}[r] 
     & *{\rule{0.1ex}{0.8ex}} \ar@{-}[r] 
     & *{\rule{0.1ex}{0.8ex}} \ar@{-}[r]
     & *{\rule{0.1ex}{0.8ex}} \ar@{-}[r] \ar@/^-1.0pc/@{-}[ll]
     & *{\rule{0.1ex}{0.8ex}} \ar@{-}[r] \ar@/^-1.5pc/@{-}[lll]
     & *{\rule{0.1ex}{0.8ex}} \ar@{-}[r] \ar@/^-2.0pc/@{-}[llll]
     & *{\rule{0.1ex}{0.8ex}} \ar@{-}[r] \ar@/^-2.5pc/@{-}[lllll]
     & *{}\ar@{--}[r]
     & *{}
                    }
}
\end{equation}
This gives an example of a maximal $1$-orthogonal subcategory of $\sD$
which is not functorially finite; see Example
\ref{exa:non-functorially_finite}.

If $\subcat$ is a cluster tilting subcategory of $\sD$, then we will
call the collection $\cluster$ of indecomposable objects of $\subcat$
a cluster.  Since $\subcat = \add \cluster$, the subcategory and the
corresponding cluster contain the same information.  Our third main
result is the following.

{\bf Theorem C.}
{\em 
The clusters form a cluster structure in $\sD$.
}

The notion of a cluster structure was introduced in \cite{BIRS} and we
have reproduced it in Section \ref{sec:cluster}.  Some of the salient
features are that if $\cluster$ is a cluster and $a$ is an
indecomposable object in $\cluster$, then $a$ can be replaced with a
unique other indecomposable object $a^*$ of $\sD$ such that a new
cluster $\cluster^*$ results, and that passing from the
Auslander-Reiten (AR) quiver of $\add \cluster$ to the AR quiver of
$\add \cluster^*$ is given by Fomin-Zelevinsky mutation at $a$ in the
sense of \cite[sec.\ 8]{FZ2}.

There are several viewpoints on the results of this paper.

(1) As mentioned, $\sD$ behaves like a cluster category of type
$A_{\infty}$.

(2) Theorems A and B show that $\sD$ can be viewed as a
categorification of triangulations of the $\infty$-gon.  Such
triangulations have not, to our knowledge, been studied elsewhere, but
they seem to be interesting combinatorial objects.

(3) Theorem C shows that $\sD$ provides a cluster tilting theory for
the abelian ca\-te\-go\-ri\-es of the form $\sD/\subcat$ where
$\subcat$ is a cluster tilting subcategory of $\sD$.  Namely, we have
$\sD/\subcat \simeq \mod\,\subcat$ by \cite[sec.\ 2]{KellerReiten2}
and \cite[cor.\ 4.4]{KoenigZhu}, and Theorem C says that the AR
quivers of $\subcat = \add \cluster$ and $\subcat^* =
\add \cluster^*$ are related by Fomin-Zelevinsky mutation at $a$, so
passing from $\sD/\subcat \simeq \mod\,\subcat$ to $\sD/\subcat^*
\simeq \mod\,\subcat^*$ can be viewed as `cluster tilting at $a$'.
Some of the categories $\mod\,\subcat$ are hereditary categories of
the form $\rep \Gamma$ where $\Gamma$ is an infinite quiver; see
Example \ref{exa:Gamma}.  Such categories were investigated by Reiten
and Van den Bergh in \cite{RVdB} and form an important branch in the
taxonomy of hereditary categories.

Other aspects of the category $\sD$ have been studied in the
literature: It is e\-qui\-va\-lent to the category $\cC_{\cH}$ which
appeared in \cite[sec.\ 2.1]{KellerReiten2} where a cluster tilting
sub\-ca\-te\-go\-ry was also shown, its Hall algebra was computed in
\cite{KYZ}, and some relations with algebraic topology were
investigated in \cite{K}.

Let us remark that we will actually obtain $\sD$ by a different recipe
from the one mentioned above.  Namely, we will use $\sD = \Df(k[T])$
where $k[T]$ is viewed as a DG algebra with $T$ placed in homological
degree $1$ and zero differential.  This gives a category which is
triangulated equivalent to the one above by \cite[sec.\ 8]{J}, but
some computations become easier.

The paper is organised as follows.  Section \ref{sec:D} gives basic
information on the category $\sD$.  Section \ref{sec:Hom} investigates
the morphisms of $\sD$.  Section \ref{sec:arcs} gives the information
we need on triangulations of the $\infty$-gon.  Section
\ref{sec:cluster_tilting} proves Theorems A and B, and Section
\ref{sec:cluster} proves Theorem C.

Section \ref{sec:questions} presents some questions; for instance, is
it possible to define a `cluster algebra of type $A_{\infty}$'?

\begin{Notation}
The set of morphisms in $\sD$ from $x$ to $y$ is denoted $\sD(x,y)$.

We will join a common abuse of terminology by saying `indecomposable
object' when we mean `isomorphism class of indecomposable objects',
and by viewing two subcategories of $\sD$ as equal if they have the
same essential closure.
\end{Notation}

\section{Basic properties of the category $\sD$}
\label{sec:D}

This section defines the category $\sD$ and recalls a few basic
properties.

\begin{Setup}
\label{set:blanket}
Throughout, $k$ is a field and $R = k[T]$ is the polynomial algebra.
We view $R$ as a DG algebra with zero differential and $T$ placed in
homological degree $1$.

Our main object of study is
\[
  \sD = \Df(R),
\]
the derived category of DG $R$-modules with finite dimensional
homology over $k$.  The suspension, Serre functor, and AR quiver of
$\sD$ will be denoted by $\Sigma$, $S$, and $Q$.
\end{Setup}

The next three remarks sum up some results on $\sD$ from \cite[section
8 in particular]{J}.

\begin{Remark}
[Basic properties]
\label{rmk:basics}
The category $\sD$ has finite dimensional $\Hom$ spaces over $k$ and
split idempotents, so it is Krull-Schmidt.  It is a $2$-Calabi-Yau
triangulated category, that is, its Serre functor is $S = \Sigma^2$
where $\Sigma$ denotes the suspension functor.  Accordingly, the AR
translation is $\tau = S\Sigma^{-1} = \Sigma$.
\end{Remark}

\begin{Remark}
[Indecomposable objects]
\label{rmk:Y}
For each integer $r \geq 0$, we have a DG $R$-module
\[
  X_r = R/(T^{r+1})
\]
which is concentrated in homological degrees from $0$ to $r$.  The
indecomposable objects of $\sD$ are $\Sigma^j X_r$ for $j$, $r$
integers, $r \geq 0$.

There is an obvious short exact sequence of DG modules $0 \rightarrow
\Sigma^{r+1}R \rightarrow R \rightarrow X_r \rightarrow 0$ which
induces a distinguished triangle
\begin{equation}
\label{equ:a}
  \Sigma^{r+1}R
  \rightarrow R
  \rightarrow X_r
  \rightarrow \Sigma^{r+2}R
\end{equation}
in $\sD$.  Hence the DG module $X_r$ is quasi-isomorphic to the
mapping cone $C_r$ of $\Sigma^{r+1}R \rightarrow R$.

Denote by $(-)^{\natural}$ the operation of forgetting the
differential.  Then $R^{\natural}$ is a graded algebra,
$C_r^{\natural}$ is a graded $R^{\natural}$-module, and the
construction of the mapping cone gives
\[
  C_r^{\natural} = R^{\natural} \oplus \Sigma^{r+2}R^{\natural}.
\]
Denoting the generators of the two copies of $R^{\natural}$ by $e_0$
and $e_{r+2}$, the differential of $C_r$ is given by
\[
  \partial(e_0) = 0,  \;\;\;  \partial(e_{r+2}) = T^{r+1}e_0.
\]
It is easy to see that $C_r$ is a minimal semi-free resolution of
$X_r$.
\end{Remark}

\begin{Remark}
[Auslander-Reiten quiver]
\label{rmk:Q}
The AR quiver $Q$ of $\sD$ is $\BZ A_{\infty}$ and the indecomposable
objects are arranged in the quiver as follows.
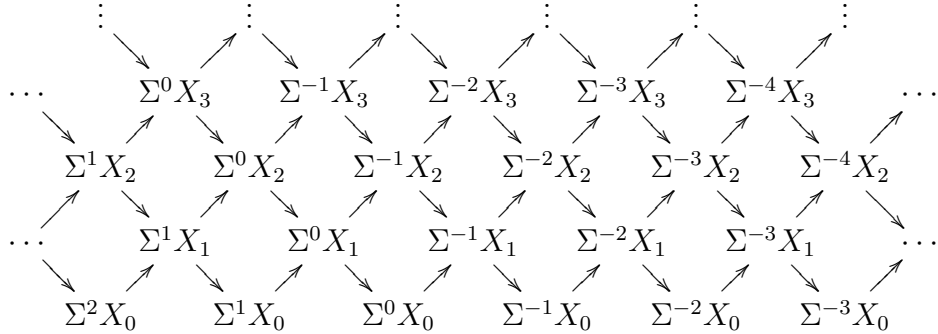
\[
  \xymatrix @-3.1pc @! {
    & \vdots \ar[dr] & & \vdots \ar[dr] & & \vdots \ar[dr] & & \vdots \ar[dr] & & \vdots \ar[dr] & & \vdots & \\
    \cdots \ar[dr]& & \Sigma^0 X_3 \ar[ur] \ar[dr] & & \Sigma^{-1} X_3 \ar[ur] \ar[dr] & & \Sigma^{-2} X_3 \ar[ur] \ar[dr] & & \Sigma^{-3}X_3 \ar[ur] \ar[dr] & & \Sigma^{-4}X_3 \ar[ur] \ar[dr] & & \cdots \\
    & \Sigma^1 X_2 \ar[ur] \ar[dr] & & \Sigma^0 X_2 \ar[ur] \ar[dr] & & \Sigma^{-1} X_2 \ar[ur] \ar[dr] & & \Sigma^{-2}X_2 \ar[ur] \ar[dr] & & \Sigma^{-3}X_2 \ar[ur] \ar[dr] & & \Sigma^{-4}X_2 \ar[ur] \ar[dr] & \\
    \cdots \ar[ur]\ar[dr]& & \Sigma^1 X_1 \ar[ur] \ar[dr] & & \Sigma^0 X_1 \ar[ur] \ar[dr] & & \Sigma^{-1} X_1 \ar[ur] \ar[dr] & & \Sigma^{-2}X_1 \ar[ur] \ar[dr] & & \Sigma^{-3}X_1 \ar[ur] \ar[dr] & & \cdots\\
    & \Sigma^2 X_0 \ar[ur] & & \Sigma^1 X_0 \ar[ur] & & \Sigma^0 X_0 \ar[ur] & & \Sigma^{-1}X_0 \ar[ur] & & \Sigma^{-2}X_0 \ar[ur] & & \Sigma^{-3}X_0 \ar[ur] & \\
               }
\]
We will use the following standard coordinate system on $Q$.
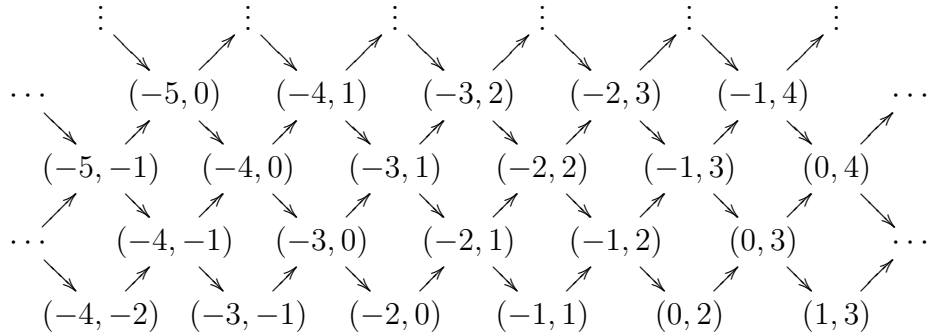
\[
  \xymatrix @-4.0pc @! {
    & \vdots \ar[dr] & & \vdots \ar[dr] & & \vdots \ar[dr] & & \vdots \ar[dr] & & \vdots \ar[dr] & & \vdots & \\
    \cdots \ar[dr]& & (-5,0) \ar[ur] \ar[dr] & & (-4,1) \ar[ur] \ar[dr] & & (-3,2) \ar[ur] \ar[dr] & & (-2,3) \ar[ur] \ar[dr] & & (-1,4) \ar[ur] \ar[dr] & & \cdots \\
    & (-5,-1) \ar[ur] \ar[dr] & & (-4,0) \ar[ur] \ar[dr] & & (-3,1) \ar[ur] \ar[dr] & & (-2,2) \ar[ur] \ar[dr] & & (-1,3) \ar[ur] \ar[dr] & & (0,4) \ar[ur] \ar[dr] & \\
    \cdots \ar[ur]\ar[dr]& & (-4,-1) \ar[ur] \ar[dr] & & (-3,0) \ar[ur] \ar[dr] & & (-2,1) \ar[ur] \ar[dr] & & (-1,2) \ar[ur] \ar[dr] & & (0,3) \ar[ur] \ar[dr] & & \cdots\\
    & (-4,-2) \ar[ur] & & (-3,-1) \ar[ur] & & (-2,0) \ar[ur] & & (-1,1) \ar[ur] & & (0,2) \ar[ur] & & (1,3) \ar[ur] & \\
               }
\]
Accordingly, coordinate pairs and indecomposable objects will be
related by
\[
  (m,n) = \Sigma^{-n}X_{n-m-2}.
\]
Note that in terms of coordinates, the actions of $\Sigma = \tau$ and
$S = \Sigma^2$ on objects are
\[
  \Sigma(m,n) = (m-1,n-1),  \;\;  S(m,n) = (m-2,n-2).
\]
\end{Remark}

\section{Morphisms in the category $\sD$}
\label{sec:Hom}

This technical section provides detailed information on the morphisms
of the category $\sD$.

\begin{Definition}
\label{def:H}
Let $x = (i,j)$ be a vertex of the AR quiver $Q$ of $\sD$.  We define
(unbounded) subsets $H^-(x)$ and $H^+(x)$ of vertices of $Q$ which can
be sketched as follows.
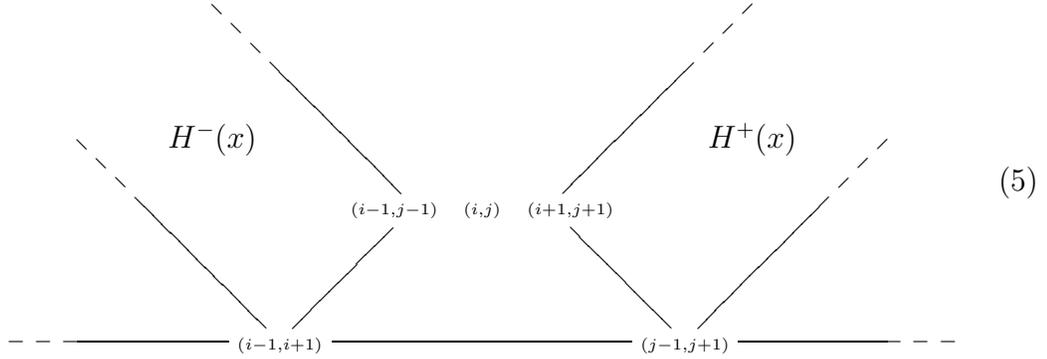
\begin{equation}
\label{fig:H}
\vcenter{
  \xymatrix @-4.5pc @! {
    &&&*{} &&&&&&&& *{}&& \\
    &&&& *{} \ar@{--}[ul] & & & & & & *{} \ar@{--}[ur] \\
    &*{}&& H^-(x) & & & & & & & & H^+(x) && *{}\\
    &&*{}\ar@{--}[ul]& & & & {\scriptscriptstyle (i-1,j-1)\hspace{3ex}} \ar@{-}[ddll] \ar@{-}[uull] & {\scriptscriptstyle (i,j)} & {\hspace{3ex}\scriptscriptstyle (i+1,j+1)} \ar@{-}[ddrr] \ar@{-}[uurr]& & &&*{}\ar@{--}[ur]&\\ 
    && \\
    *{}\ar@{--}[r]&*{} \ar@{-}[rrr] && & {\scriptscriptstyle (i-1,i+1)} \ar@{-}[uull]\ar@{-}[rrrrrr]& & & & & & {\scriptscriptstyle (j-1,j+1)} \ar@{-}[uurr]\ar@{-}[rrr]&&&*{}\ar@{--}[r]&*{}\\
           }
}
\end{equation}
The subsets include the edges.  In a more rigorous vein, we have
\begin{align*}
  H^-(x) & = \{\, (m,n) \:|\: m \leq i-1, \;\: i+1 \leq n \leq j-1 \,\}, \\
  H^+(x) & = \{\, (m,n) \:|\: i+1 \leq m \leq j-1, \;\: j+1 \leq n \,\}.
\end{align*}
We write $H(x) = H^-(x) \cup H^+(x)$.
\end{Definition}

%The following Proposition says that the indecomposable objects $y$ to
%which $x$ has non-zero morphisms form two regions as follows.
%\[
%  \xymatrix @-4.0pc @! {
%    &&&*{} &&&&&&&& *{}&& \\
%    &&&& *{} \ar@{--}[ul] & & & & & & *{} \ar@{--}[ur] \\
%    &*{}&& H^-(\Sigma x) & & & & & & & & H^+(\Sigma x) && *{}\\
%    &&*{}\ar@{--}[ul]& & & & \Sigma^2 x \ar@{-}[ddll] \ar@{-}[uull] & \Sigma x & x \ar@{-}[ddrr] \ar@{-}[uurr]& & &&*{}\ar@{--}[ur]&\\ 
%    && \\
%    *{}\ar@{--}[r]&*{} \ar@{-}[rrr] && & *{} \ar@{-}[uull]\ar@{-}[rrrrrr]& & & & & & *{} \ar@{-}[uurr]\ar@{-}[rrr]&&&*{}\ar@{--}[r]&*{}\\
%           }
%\]

The following Proposition says that an indecomposable object $x$ has
non-zero morphisms to the indecomposable objects $y$ in two regions
like the ones in figure \eqref{fig:H}.  The object $x$ is at the
leftmost vertex of the right hand region.

\begin{Proposition}
\label{pro:Hom}
Let $x$ and $y$ be indecomposable objects of $\sD$.  Then
\[
  \sD(x,y)
  = \left\{
      \begin{array}{cl}
        k & \mbox{for } y \in H(\Sigma x), \\
        0 & \mbox{otherwise.}
      \end{array}
    \right.
\]
\end{Proposition}

\begin{proof}
Using a power of $\Sigma$, we can shift $x$ and $y$ horizontally on
the quiver without loss of generality, and so we can assume
\[
  x = (-r-2,0) = X_r.
\]
We will write
\[
  y = (m,n) = \Sigma^{-n}X_{n-m-2}.
\]

Taking $\Hom$ of the distinguished triangle \eqref{equ:a} into the
object $\Sigma^{-n}X_{n-m-2}$ gives a long exact sequence containing
\begin{align*}
  & \sD(\Sigma^{r+2}R,\Sigma^{-n}X_{n-m-2})
  \rightarrow \sD(X_r,\Sigma^{-n}X_{n-m-2}) \\
  & \hspace{20ex} \rightarrow \sD(R,\Sigma^{-n}X_{n-m-2})
  \rightarrow \sD(\Sigma^{r+1}R,\Sigma^{-n}X_{n-m-2})
\end{align*}
which is
\begin{align}
\label{equ:b}
  & \H_{n+r+2}(X_{n-m-2})
  \rightarrow \sD(x,y) \\
\nonumber
  & \hspace{15ex} \rightarrow \H_n(X_{n-m-2})
  \rightarrow \H_{n+r+1}(X_{n-m-2}).
\end{align}

Consider the sketch of the AR quiver below.  It is cumbersome, but
elementary, to verify from \eqref{equ:b} that $\sD(x,y)$ is $k$ when
$y$ is in the region $\xymatrix{*+[o][F]{2}}$ (which includes the
edges and is equal to $H^+(\Sigma x)$).  Also, $\sD(x,y)$ is $0$ when
$y$ is in region $\xymatrix{*+[o][F]{1}}$ (which does not include the
diagonal edge) or region $\xymatrix{*+[o][F]{3}}$ (which includes the
dotted edge, but not the other one).
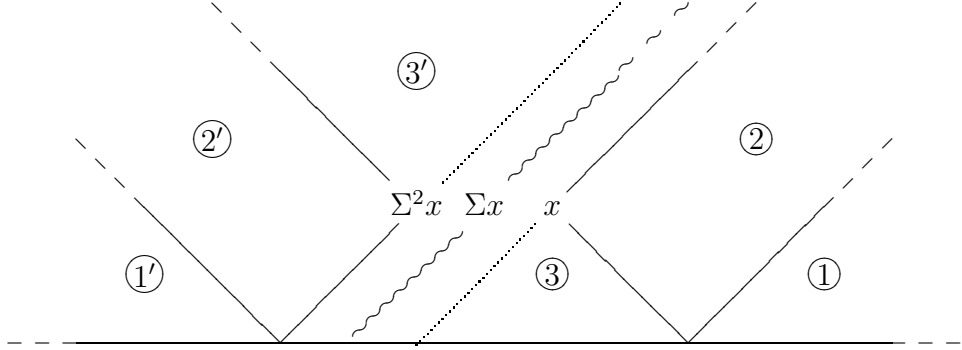
\[
  \xymatrix @-2.1pc @! {
    &&&*{} &&&&&& *{} & *{} & *{}&& \\
    && & & *{}\ar@{--}[ul] & & *+[o][F]{3^{\prime}} & & & *{}\ar@{~~}[ur]& *{} \ar@{--}[ur]\\
    &*{}&& *+[o][F]{2^{\prime}} & & & & & & & & *+[o][F]{2} && *{}\\
    &&*{}\ar@{--}[ul]& & & & \Sigma^2 x \ar@{-}[ddll] \ar@{-}[uull] \ar@{.}[uuurrr]& \Sigma x \ar@{~}[uurr]& x \ar@{-}[ddrr] \ar@{-}[uurr]& & &&*{}\ar@{--}[ur]&\\ 
    && *+[o][F]{1^{\prime}} & & & & & & *+[o][F]{3} & & & & *+[o][F]{1} \\
    *{}\ar@{--}[r]&*{} \ar@{-}[rrr] && & *{} \ar@{-}[uull]\ar@{-}[rr]& *{} \ar@{~}[uurr]& *{} \ar@{.}[uurr]\ar@{-}[rrrr]& & & & *{} \ar@{-}[uurr]\ar@{-}[rrr]&&&*{}\ar@{--}[r]&*{}\\
           }
\]
Serre duality says $\sD(a,b) \cong \sD(b,\Sigma^2 a)^{\vee}$.  By
applying this to the previous results, we get that $\sD(x,y)$ is $k$
when $y$ is in the region $\xymatrix{*+[o][F]{2^{\prime}}}$ (which
includes the edges and is equal to $H^-(\Sigma x)$).  Also, $\sD(x,y)$
is $0$ when $y$ is in region $\xymatrix{*+[o][F]{1^{\prime}}}$ (which
does not include the diagonal edge) or region
$\xymatrix{*+[o][F]{3^{\prime}}}$ (which includes the dotted edge but
not the other one).

To complete the proof, we must show that $\sD(x,y)$ is $0$ when
$y$ is on the wavy line through $\Sigma x$.  The vertices on this
line have the form $(-r-3,-r-1+t)$ for $t \geq 0$, that is, they are
the objects $\Sigma^{r+1-t}X_t$ for $t \geq 0$, and we must show that
a morphism $X_r \rightarrow \Sigma^{r+1-t}X_t$ in $\sD$ is $0$.  Such
a morphism is a homotopy class of morphisms of DG modules
\[
  \gamma: C_r \rightarrow \Sigma^{r+1-t}X_t
\]
where $C_r$ is the minimal semi-free resolution of $X_r$ from Remark
\ref{rmk:Y}.

Recall that $C_r$ has generators $e_0$ and $e_{r+2}$ in homological
degrees $0$ and $r+2$.  The DG module $\Sigma^{r+1-t}X_t$ is
concentrated in homological degrees from $r+1-t$ to $r+1$.

If $r+1-t > 0$, then $\Sigma^{r+1-t}X_t$ is $0$ in each degree where
$C_r$ has a generator, so $\gamma = 0$ is clear.

If $r+1-t \leq 0$, then $\Sigma^{r+1-t}X_t$ is $0$ in degree $r+2$,
but it is $k$ in degree $0$.  Potentially, $\gamma(e_0)$ could be
non-zero.  However, 
\[
  T^{r+1}\gamma(e_0) = \gamma(T^{r+1}e_0) = \gamma\partial(e_{r+2})
  = \partial\gamma(e_{r+2})
  \stackrel{\rm (a)}{=} \partial(0) = 0
\]
where (a) is because $\Sigma^{r+1-t}X_t$ is $0$ in degree $r+2$.
Since $\Sigma^{r+1-t}X_t$ is equal to $k$ in degree $r+1$, it follows
that $\gamma(e_0) = 0$ so $\gamma = 0$.
\end{proof}

\begin{Corollary}
\label{cor:to_and_from}
Let $x$ and $y$ be indecomposable objects of $\sD$.  The following are
equivalent.
\begin{enumerate}
  
  \item  $\sD(x,y) \neq 0$.

\smallskip

  \item  $\sD(x,y) = k$.

\smallskip

  \item  $y \in H(\Sigma x)$.

\smallskip

  \item  $x \in H(\Sigma^{-1}y)$.

\end{enumerate}
\end{Corollary}

\begin{proof}
(i), (ii), and (iii) are equivalent by Proposition \ref{pro:Hom}.
Using Serre duality, (i) is equivalent to $\sD(y,\Sigma^2 x) \neq 0$.
Using (iii), this is equivalent to $\Sigma^2 x \in H(\Sigma y)$, that
is $x \in H(\Sigma^{-1}y)$, and this is (iv).
\end{proof}

\begin{Remark}
[Forward morphisms]
\label{rmk:forward}
Proposition \ref{pro:Hom} and Corollary \ref{cor:to_and_from} show
that there are two distinct types of non-zero morphisms going from $x$
to indecomposable objects of $\sD$: Those going to objects in
$H^+(\Sigma x)$ will be called {\em forward morphisms}, and those
going to objects in $H^-(\Sigma x)$ will be called {\em backward
morphisms}.  The backward morphisms cannot be seen in the AR quiver;
they are in the infinite radical of $\sD$.

The forward morphisms have an easy model: Up to multiplication by a
non-zero scalar, they are induced by certain canonical morphisms of DG
modules.  Namely, consider again the case $x = (-r-2,0) = X_r$.  Then
$x$ is a DG module which is concentrated in homological degrees from
$0$ to $r$.  Let $y = (m,n) = \Sigma^{-n}X_{n-m-2}$ be in the region
$H^+(\Sigma x)$ whence $-r-2 \leq m \leq -2$ and $n \geq 0$.  Then $y$
is a DG module which is concentrated in homological degrees from $-n$
to $-m-2$, and we have $-n \leq 0$ and $0 \leq -m-2 \leq r$.  We can
sketch the non-zero parts of the DG modules $x = X_r$ and $y =
\Sigma^{-n}X_{n-m-2}$ as follows, where the numbers at the top are
homological degrees and where each horizontal line indicates the
degrees where a module has non-zero components.
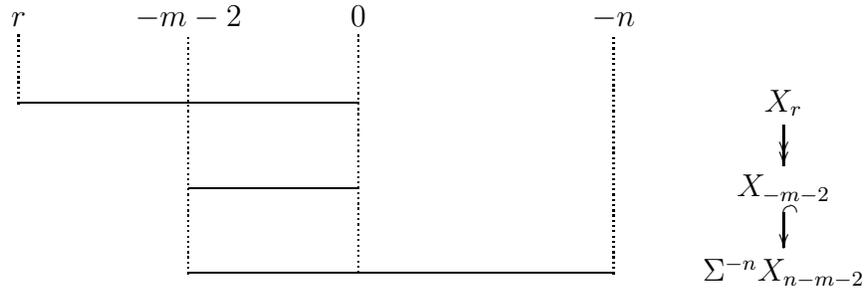
\[
  \xymatrix @-5pc @! {
    r \ar@{.}[d]   & & -m-2 \ar@{.}[ddd]  & & 0 \ar@{.}[ddd] & & & -n \ar@{.}[ddd] & & \\
    *{} \ar@{-}[rrrr] & &                   & & *{} & & &     & & X_r \ar@{->>}[d]\\
                     & & *{} \ar@{-}[rr]    & & *{} & & &     & & X_{-m-2} \ar@{^{(}->}[d]\\
                     & & *{} \ar@{-}[rrrrr] & &     & & & *{} & & \Sigma^{-n}X_{n-m-2}\\
           }
\]
We have included the DG module $X_{-m-2}$ in the sketch, and it is
clear that there is a surjective and an injective morphism of DG
modules as indicated.  Their composition is a canonical morphism of DG
modules which induces a forward morphism $x \rightarrow y$ in $\sD$.

Observe that the canonical morphism of DG modules induces a (non-zero)
forward morphism $x \rightarrow y$ in $\sD$ if and only if there is a
homological degree where both $x$ and $y$ have a non-zero component.
Indeed, it is easy to check that this happens if and only if $y$ is in
$H^+(\Sigma x)$.
\end{Remark}

\begin{Lemma}
\label{lem:xyz1}
Let $x$, $y$, and $z$ be indecomposable objects of $\sD$ such that
$y,z \in H^+(\Sigma x)$ and $z \in H^+(\Sigma y)$, for instance as in the
following sketch.
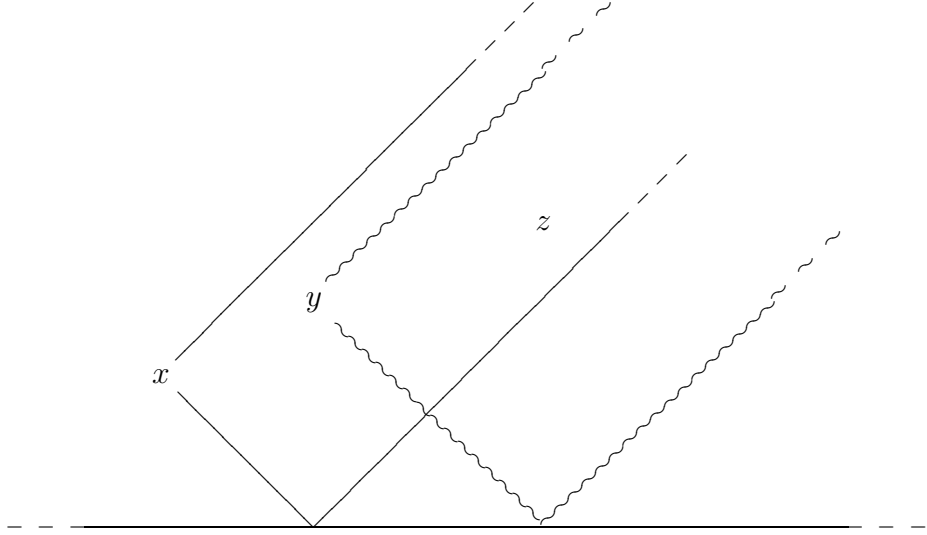
\[
  \xymatrix @-0.8pc @! {
    & & & & & & & & \\
    & & & & & & *{} \ar@{--}[ur] & *{} \ar@{~~}[ur] & *{} \\
    & & & & & & & & & \\
    & & & & & & & z &  *{}\ar@{--}[ur] & & & & \\
    & & & & y \ar@{~}[dddrrr] \ar@{~}[uuurrr] & & & & & *{} & *{}\ar@{~~}[ur]\\
    *{} & & x \ar@{-}[ddrr] \ar@{-}[uuuurrrr]& & &&&&&\\ 
    & & & & & & & *{} \\
    *{} \ar@{--}[r] & *{} \ar@{-}[rrr] & & & *{} \ar@{-}[uuuurrrr]\ar@{-}[rrrrrrr]&& &*{}\ar@{~}[uuurrr]&&&&*{}\ar@{--}[r]&*{}\\
           }
\]
\begin{enumerate}

  \item  The composition of non-zero morphisms $x \rightarrow y$ and $y
  \rightarrow z$ is non-zero.

\smallskip

  \item  Let $y \stackrel{f}{\rightarrow} z$ be a non-zero morphism.
Then each morphism $x \rightarrow z$ factors as $x \rightarrow y
\stackrel{f}{\rightarrow} z$.

\end{enumerate}
\end{Lemma}

\begin{proof}
(i).  Since we have $y \in H^+(\Sigma x)$ and $z \in H^+(\Sigma y)$,
the non-zero morphisms $x \rightarrow y$ and $y \rightarrow z$ in
$\sD$ are forward morphisms.  By Remark \ref{rmk:forward}, up to
multiplication by non-zero scalars which can be ignored, they are
induced by canonical morphisms of DG modules which we can indicate as
follows.
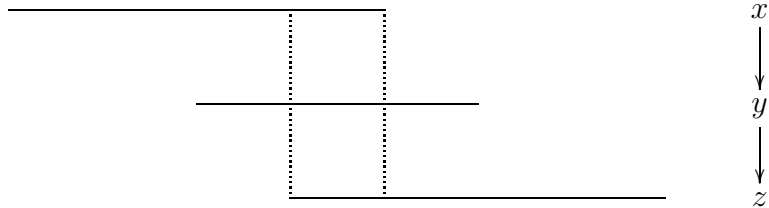
\[
  \xymatrix @-0.25pc @! {
    *{} \ar@{-}[rrrr] & & & *{} \ar@{.}[dd] & *{} \ar@{.}[dd] & & &     & x \ar[d]\\
                     & & *{} \ar@{-}[rrr]    & & *{} & *{} & &     & y \ar[d]\\
                     & & & *{} \ar@{-}[rrrr] & *{} & & & *{} & z \\
           }
\]
It is clear that these compose to a canonical morphism $x \rightarrow
z$.

Since we have $z \in H^+(\Sigma x)$, Remark \ref{rmk:forward} gives
that there is a homological degree where both the DG modules $x$ and
$z$ have a non-zero component, and hence the canonical morphism $x
\rightarrow z$ induces a (non-zero) forward morphism $x \rightarrow z$
in $\sD$ as desired.

(ii).  We must show that $\sD(x,f): \sD(x,y) \rightarrow \sD(x,z)$ is
surjective.  Since each non-zero $\Hom$ set is isomorphic to $k$, it
is enough to see that $\sD(x,f)$ is non-zero, and this follows from
part (i) because it sends $x \rightarrow y$ to the composition $x
\rightarrow y \stackrel{f}{\rightarrow} z$.
\end{proof}

\begin{Lemma}
\label{lem:S}
Let $x$, $y$, and $z$ be indecomposable objects of $\sD$.
\begin{enumerate}

  \item  $y \in H^+(\Sigma x) \Leftrightarrow Sx \in H^-(\Sigma y)$.

  \item  $z \in H^-(\Sigma x) \Leftrightarrow Sx \in H^+(\Sigma z)$.

\end{enumerate}
\end{Lemma}

\begin{proof}
(i) Suppose $y \in H^+(\Sigma x)$; then there is a non-zero morphism
$x \rightarrow y$ in $\sD$ by Corollary \ref{cor:to_and_from}.  The
Serre duality isomorphism $\sD(x,y) \cong \sD(y,Sx)^{\vee}$ implies
that there is a non-zero morphism $y \rightarrow Sx$ so we have $Sx
\in H(\Sigma y)$.  To establish the implication $\Rightarrow$, it
remains to see that $Sx$ is in $H^-(\Sigma y)$, not $H^+(\Sigma y)$.

However, if $x = (i,j)$ then the shape of the region $H^+(\Sigma x)$
implies $y = (i+p,j+q)$ for some $p,q \geq 0$.  This again means that
the points of $H^+(\Sigma y)$ have the form
$(i+p+p^{\prime},j+q+q^{\prime})$ for some $p^{\prime}, q^{\prime}
\geq 0$, but $Sx = (i-2,j-2)$ is not of this form so we must have $Sx$
in $H^-(\Sigma y)$.

The implication $\Leftarrow$ is proved by a similar argument.

(ii)  We have
\[
  z \in H^-(\Sigma x)
  \Leftrightarrow SS^{-1}z \in H^-(\Sigma x)
  \Leftrightarrow x \in H^+(\Sigma S^{-1}z)
  \Leftrightarrow Sx \in H^+(\Sigma z)
\]
where the second biimplication is by part (i).
\end{proof}

\begin{Lemma}
\label{lem:xyz2}
Let $x$, $y$, and $z$ be indecomposable objects of $\sD$ such that
$y, z \in H^-(\Sigma x)$ and $z \in H^+(\Sigma y)$, for instance as in
the following sketch.
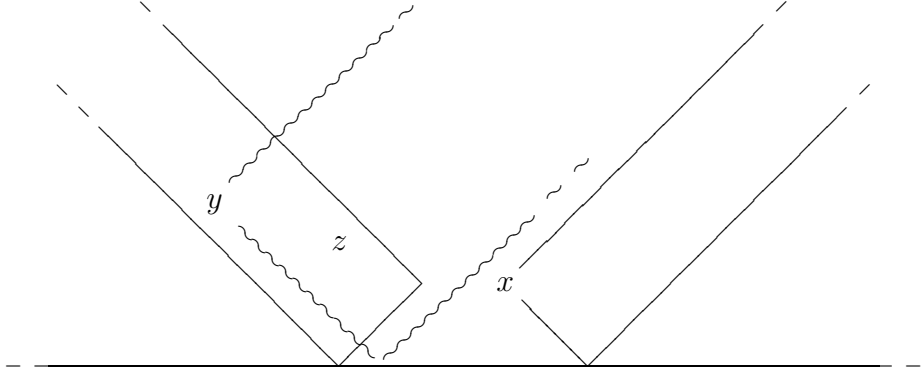
\[
  \xymatrix @-1.9pc @! {
    & & & & & & & & & & & & & & & & & & & & & & \\
    & & & & *{} \ar@{--}[ul] & & & & & *{} \ar@{~~}[ur] & & & & & & & & & *{} \ar@{--}[ur] & \\
    & & & & & & & & & & & & & & & & & & & & & \\
    & & \ar@{--}[ul] & & & & & & & & & & & & & & & & & & *{} \ar@{--}[ur] \\
    & & & & & & &*{} & *{} & & & & & & *{} & *{} & & & & \\
    & & & & & y \ar@{~}[uuuurrrr] & & & & & & & & \ar@{~~}[ur] & & & & & \\
    & & & & & & & & z & & & & & & & & & \\
    & & & & & &*{}& & & & *{} \ar@{-}[ddll] \ar@{-}[uuuuuullllll] & & x \ar@{-}[ddrr] \ar@{-}[uuuuuurrrrrr]& & &&*{}&\\ 
    & & & & & & \\
    *{} \ar@{--}[r] & *{} \ar@{-}[rrrrrrr] &&&&&& & *{} \ar@{-}[uuuuuullllll]\ar@{-}[rrrrrr]& *{} \ar@{~}[uuuurrrr]\ar@{~}[uuuullll]& & & & & *{} \ar@{-}[uuuuuurrrrrr]\ar@{-}[rrrrrrr]&&&&&&&*{}\ar@{--}[r]&*{}\\
           }
\]
Let $y \stackrel{f}{\rightarrow} z$ be a non-zero morphism.  Then each
morphism $x \rightarrow z$ factors as $x
\rightarrow y \stackrel{f}{\rightarrow} z$.
\end{Lemma}

\begin{proof}
We must show that $\sD(x,f): \sD(x,y) \rightarrow \sD(x,z)$ is
surjective.  Using Serre duality, it is the same to show that
$\sD(f,Sx): \sD(z,Sx) \rightarrow \sD(y,Sx)$ is injective.  This map
sends $z \stackrel{\zeta}{\rightarrow} Sx$ to the composition $y
\stackrel{f}{\rightarrow} z \stackrel{\zeta}{\rightarrow} Sx$.
However, we have $z \in H^-(\Sigma x)$ so Lemma \ref{lem:S}(ii) says
$Sx \in H^+(\Sigma z)$.  Hence, if $\zeta$ is non-zero then it is a
forward morphism.  So is $f$ since $z \in H^+(\Sigma y)$, and then $y
\stackrel{f}{\rightarrow} z \stackrel{\zeta}{\rightarrow} Sx$ is
non-zero by Lemma \ref{lem:xyz1}(i) since we have $Sx \in \H^+(\Sigma
y)$; this holds by Lemma \ref{lem:S}(ii) again since $y \in H^-(\Sigma
x)$.
\end{proof}

\section{Triangulations of the $\infty$-gon}
\label{sec:arcs}

This section studies triangulations of the $\infty$-gon, that is,
maximal sets of non-intersecting arcs, and their relation with the
category $\sD$.

%\[
%  \rule{0ex}{10ex}
%  \xymatrix @+0.1pc @! {
%       *{} \ar@{--}[r]
%     & *{}\ar@{-}[r]
%     & *{\rule{0.1ex}{0.8ex}} \ar@{-}[r]
%     & *{\rule{0.1ex}{0.8ex}} \ar@{-}[r]
%     & *{\rule{0.1ex}{0.8ex}} \ar@{-}[r]
%     & *{\rule{0.1ex}{0.8ex}} \ar@{-}[r]
%     & *{\rule{0.1ex}{0.8ex}} \ar@{-}[r]
%     & *{\rule{0.1ex}{0.8ex}} \ar@{-}[r]
%     & *{\rule{0.1ex}{0.8ex}} \ar@{-}[r]
%     & *{}\ar@{--}[r]
%     & *{}
%                    }
%\]

\begin{Definition}
An {\em arc} is a pair $(m,n)$ of integers with $m \leq n-2$.

The arc $(m,n)$ is said to {\em end} in each of the integers $m$ and
$n$.
%and to {\em connect} the (non-neighbouring) integers $m$ and $n$.

Two arcs $(m,n)$ and $(p,q)$ are said to {\em intersect} if we have
either $m < p < n < q$ or $p < m < q < n$.
\end{Definition}

The definition is intended to capture our geometric intuition in which
an arc is drawn as a curve between two integers on the number line as
follows.
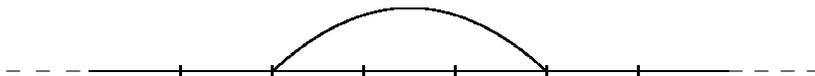
\[
  \xymatrix @-2.5pc @! {
    \rule{0ex}{6.5ex}\ar@{--}[r]&\ar@{-}[r]&*{\rule{0.1ex}{0.8ex}} \ar@{-}[r]& *{\rule{0.1ex}{0.8ex}} \ar@{-}[r] \ar@/^2pc/@{-}[rrr]& *{\rule{0.1ex}{0.8ex}} \ar@{-}[r] & *{\rule{0.1ex}{0.8ex}} \ar@{-}[r] & *{\rule{0.1ex}{0.8ex}} \ar@{-}[r] & *{\rule{0.1ex}{0.8ex}} \ar@{-}[r] &*{}\ar@{--}[r]&*{}
                    }
\]
Two arcs can be drawn as non-intersecting curves precisely if they do
not intersect in the sense of the definition, with the proviso that
curves which only meet at their end points are not viewed as
intersecting.

In an informal sense, it is reasonable to view the integers as being
the vertices of an $\infty$-gon and to view arcs as being diagonals
between vertices.  Hence a maximal set of non-intersecting arcs can be
viewed as a triangulation of the $\infty$-gon.  Some typical ways of
achieving such maximal sets are shown in sketches
\eqref{equ:leapfrog}, \eqref{equ:fountain}, and
\eqref{equ:right_and_left_fountain} in the Introduction.  The sketches
inspire the following definition.

\begin{Definition}
Let $\arcs$ be a set of arcs.  If for each integer $n$ there are only
finitely many arcs in $\arcs$ which end in $n$, then $\arcs$ is called
{\em locally finite}.

If $n$ is an integer such that $\arcs$ contains infinitely many arcs
of the form $(m,n)$, then $n$ is called a {\em left-fountain} of
$\arcs$.

If $n$ is an integer such that $\arcs$ contains infinitely many arcs
of the form $(n,p)$, then $n$ is called a {\em right-fountain} of
$\arcs$.

If $n$ is both a left- and a right-fountain of $\arcs$, then it is
called a {\em fountain}.
\end{Definition}

It turns out that if a maximal set of non-intersecting arcs has a
right-fountain then it also has a left-fountain and vice versa; we owe
this observation to Collin Bleak.  However, all we need here is the
following more modest result.

\begin{Lemma}
\label{lem:fountains}
Let $\arcs$ be a maximal set of non-intersecting arcs.  Then $\arcs$
has at most one right-fountain and at most one left-fountain.
\end{Lemma}

\begin{proof}
Suppose that $\arcs$ is not locally finite and let $m$ be an integer
where infinitely many arcs of $\arcs$ end.  So $m$ is either a right-
or a left-fountain of $\arcs$ and we can suppose the former without
loss of generality.

We must show that $m$ is the only right-fountain of $\arcs$, so let $p
\neq m$ be an integer.  If $p > m$ then we can pick $n > p$ such that
$(m,n)$ is in $\arcs$.  An arc $(p,q)$ will intersect $(m,n)$ as soon
as $q > n$.
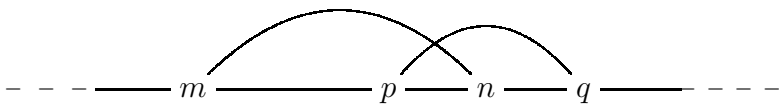
\[
  \xymatrix @-2.5pc @! {
    \rule{0ex}{6.5ex}\ar@{--}[r]&\ar@{-}[r]&{m} \ar@{-}[r]\ar@/^2.5pc/@{-}[rrr]&*{}\ar@{-}[r]& {p} \ar@{-}[r] \ar@/^2pc/@{-}[rr]& {n} \ar@{-}[r] & {q} \ar@{-}[r] & *{}\ar@{--}[r]&*{}
                    }
\]
So $\arcs$ contains only finitely many arcs of the form $(p,q)$ and
$p$ is not a right-fountain.

If $p < m$ then an arc $(p,q)$ can only be in $\arcs$ if $q \leq m$,
for if $q > m$ then it is possible to pick an arc $(m,n)$ in $\arcs$
with $n > q$, and then $(m,n)$ and $(p,q)$ intersect.
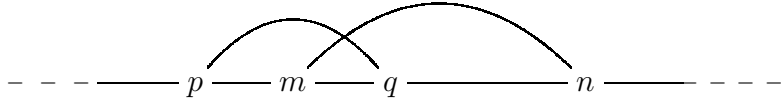
\[
  \xymatrix @-2.5pc @! {
    \rule{0ex}{6.5ex}\ar@{--}[r]&\ar@{-}[r]&{p} \ar@{-}[r]\ar@/^2pc/@{-}[rr]& {m} \ar@{-}[r] \ar@/^2.5pc/@{-}[rrr]& {q} \ar@{-}[r] & *{} \ar@{-}[r] & {n} \ar@{-}[r] & *{}\ar@{--}[r]&*{}
                    }
\]
Again $\arcs$ contains only finitely many arcs of the form $(p,q)$ and
$p$ is not a right-fountain.
\end{proof}

\begin{Remark}
\label{rmk:bijection}
An ordered pair of integers $(m,n)$ with $m \leq n-2$ can be viewed as
an arc.  Using the coordinate system of Remark \ref{rmk:Q}, it can
also be viewed as a vertex of the AR quiver $Q$ of $\sD$, that is, an
indecomposable object of $\sD$.

So there is a bijection between arcs and indecomposable objects of
$\sD$.

This induces a bijection between sets of arcs and sets of
indecomposable objects of $\sD$.  But such sets correspond bijectively
to subcategories of $\sD$ which are closed under direct sums and
direct summands, the bijection being given by $\cluster \mapsto
\add \cluster$.

So there is a bijection between sets of arcs and subcategories of
$\sD$ which are closed under direct sums and direct summands.
\end{Remark}

It is easy to check that the bijection plays together with the regions
$H(x)$ as follows.

\begin{Lemma}
\label{lem:intersect}
Let $x$ and $y$ be indecomposable objects of $\sD$.  The following
conditions are equivalent. 
\begin{enumerate}

  \item  $x \in H(y)$.

\smallskip

  \item  $y \in H(x)$.

\smallskip

  \item  The arcs corresponding to $x$ and $y$ intersect.

\end{enumerate}
\end{Lemma}

%Namely, if $y$ is the vertex $(m,n)$ then the arc corresponding to $y$
%connects the integers $m$ and $n$.  The vertices in $H^-(y)$
%correspond to arcs connecting an integer smaller than $m$ to an
%integer between $m$ and $n$.  The vertices in $H^+(y)$ correspond to
%arcs connecting an integer greater than $n$ to an integer between $m$
%and $n$.

The following is an immediate consequence.

\begin{Lemma}
\label{lem:Hom}
Let $x$ and $y$ be indecomposable objects of $\sD$.  Then
\[
  \sD(x,y) \neq 0
  \; \Leftrightarrow \;
  \mbox{the arcs corresponding to $x$ and $\Sigma^{-1}y$ intersect}.
\]
\end{Lemma}

\begin{proof}
By Corollary \ref{cor:to_and_from}, we have $\sD(x,y) \neq 0$ if
and only if $x \in H(\Sigma^{-1}y)$, and by Lemma \ref{lem:intersect},
this is the same as for the arcs corresponding to $x$ and
$\Sigma^{-1}y$ to intersect.
\end{proof}

\section{Cluster tilting subcategories of $\sD$}
\label{sec:cluster_tilting}

This section proves Theorems A and B from the Introduction; see
Theorems \ref{thm:maximal_1-orthogonal} and \ref{thm:cluster_tilting}. 

The distinction between maximal $1$-orthogonal and cluster tilting
subcategories in the following definition is not standard in the
literature, but Theorem B means that it is useful for this paper. 

\begin{Definition}
Let $\subcat$ be a subcategory of $\sD$.  We write 
\begin{align*}
  \subcat^{\perp} 
  & = \{\, d \in \sD \,|\, \sD(a,d) = 0 \;\mbox{for each}\; a \in \subcat \,\}, \\
  {}^{\perp}\subcat 
  & = \{\, d \in \sD \,|\, \sD(d,a) = 0 \;\mbox{for each}\; a \in \subcat \,\}.
\end{align*}

A subcategory $\subcat$ is called {\em maximal $1$-orthogonal} if it
satisfies $\subcat = (\Sigma^{-1}\subcat)^{\perp}$ and $\subcat =
{}^{\perp}(\Sigma \subcat)$.  (In fact, either equality implies the
other because $\sD$ is a $2$-Calabi-Yau category.)

A subcategory $\subcat$ is called {\em cluster tilting} if it is
maximal $1$-orthogonal and functorially finite.
\end{Definition}

\begin{Remark}
\label{rmk:forbidden}
Let $\subcat$ be a subcategory of $\sD$ which is closed under direct
sums and direct summands.  The inclusion
\begin{equation}
\label{equ:inclusion}
  \subcat \subseteq (\Sigma^{-1}\subcat)^{\perp}
\end{equation}
holds precisely if the presence of an indecomposable object $a$ in
$\subcat$ forbids an indecomposable object $x$ from being in $\subcat$
when there is a non-zero morphism $\Sigma^{-1}a \rightarrow x$.

It hence follows from Corollary \ref{cor:to_and_from} that the
inclusion \eqref{equ:inclusion} is equivalent to the following
condition:  If $a$ is in $\subcat$ then the indecomposable objects in
$H(\Sigma\Sigma^{-1}a) = H(a)$ are forbidden from $\subcat$.

We therefore sometimes refer to the $H(a)$ as {\em forbidden regions}.
Note that, in particular, a maximal $1$-orthogonal subcategory of
$\sD$ satisfies \eqref{equ:inclusion}.
\end{Remark}

\begin{Theorem}
\label{thm:maximal_1-orthogonal}
Let $\subcat$ be a subcategory of $\sD$ which is closed under direct
sums and direct summands.  Let $\arcs$ be the corresponding set of
arcs under the bijection of Remark \ref{rmk:bijection}.

Then $\subcat$ is a maximal $1$-orthogonal subcategory of $\sD$ if and
only if $\arcs$ is a maximal set of non-intersecting arcs.
\end{Theorem}

\begin{proof}
By Remark \ref{rmk:forbidden}, the inclusion \eqref{equ:inclusion} is
equivalent to the condition that if $a$ is in $\subcat$ then the
objects in $H(a)$ are forbidden from $\subcat$.

An indecomposable object $a$ corresponds to an arc $\fa$, and by Lemma
\ref{lem:intersect} the indecomposable objects in $H(a)$ correspond
precisely to arcs intersecting $\fa$.  So the subcategory $\subcat$
satisfies \eqref{equ:inclusion} if and only if it corresponds to a set
of non-intersecting arcs.

It follows that subcategories $\subcat$ maximal among the ones satisfying
\eqref{equ:inclusion} correspond to maximal sets of non-intersecting
arcs.  But it is not hard to check that such maximal subcategories are
precisely the ones with $\subcat = (\Sigma^{-1}\subcat)^{\perp}$, and
these are the maximal $1$-orthogonal subcategories of $\sD$.
\end{proof}

\begin{Theorem}
\label{thm:cluster_tilting}
Let $\subcat$ be a maximal $1$-orthogonal subcategory of $\sD$.  Let
$\arcs$ be the corresponding maximal set of non-intersecting arcs
under the bijection of Remark \ref{rmk:bijection}.

Then $\subcat$ is functorially finite (that is, $\subcat$ is a cluster
tilting subcategory of $\sD$) if and only if $\arcs$ is {\rm (i)}
locally finite, or {\rm (ii)} has a fountain.
\end{Theorem}

\begin{proof}
We remind the reader of Corollary \ref{cor:to_and_from} on the
relation between existence of non-zero morphisms in $\sD$ and
membership of the regions $H^+$ and $H^-$.  This will be used
repeatedly in the proof.

We must show that $\subcat$ is functorially finite if and only if
$\arcs$ satisfies condition (i) or (ii) in the theorem.  By
\cite[lem.\ 3.2.3]{KoenigZhu} and its dual, it is enough to show that
$\subcat$ is precovering or preenveloping if and only if $\arcs$
satisfies (i) or (ii).

Suppose that (i) holds.  Then it follows easily from Lemma
\ref{lem:Hom} that for each in\-de\-com\-po\-sa\-ble object $x$ of
$\sD$, only finitely many indecomposable objects of $\subcat$ have
non-zero morphisms to $x$, and this implies that $\subcat$ is
precovering.

Suppose that (i) does not hold; that is, $\arcs$ has a right- or a
left-fountain.  Without loss of generality we can suppose that $\arcs$
has a right-fountain which by Lemma \ref{lem:fountains} is the only
right-fountain of $\arcs$.  We must show that $\subcat$ is precovering
if and only if the right-fountain is also a left-fountain.

Suppose first that $\subcat$ is precovering.  The right-fountain of
$\arcs$ is an integer $n$ for which there are infinitely many arcs of
the form $(n,p)$ in $\arcs$.  These arcs give a collection $P$ of
indecomposable objects in $\subcat$ which sit on a diagonal half line
$r$ in the AR quiver $Q$ of $\sD$.  The following sketch of $Q$ shows
$r$ along with some of the indecomposable objects $a$ of $P$ and, in
dotted lines, their regions $H(\Sigma a)$.
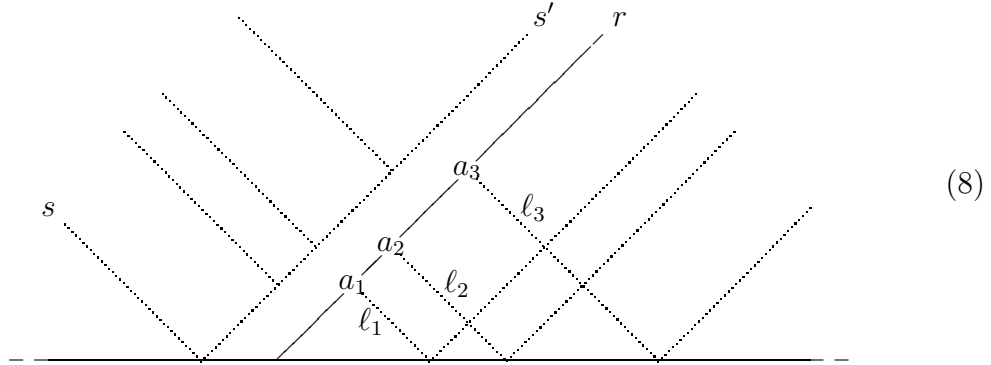
\begin{equation}
\label{equ:fig1}
\vcenter{
  \xymatrix @-3.3pc @! {
    &&&&&&*{}&&&&&&&&{s^{\prime}}&&{r} \\
    &&&&&&&&&&&&&*{}&&*{}\ar@{--}[ur] \\
    &&&&*{}&&&&&&&&*{}&&*{}&&&&*{}\\
    &&&*{}&&&&&&&&&&*{}&&&&&*{}&*{}\\
    &&&&*{} &&&&&*{}& *{} \ar@{.}[uuuullll]& *{} & *{a_3} \ar@{-}[uuurrr]\ar@{.}[dddddrrrrr] &&&&*{}&&*{} \\
    &{s}&& & & *{} & & & & & & *{} &&&{\!\!\!\! \ell_3}&*{}&*{}&*{}&&&&*{}\\
    &&*{}&& & & & & *{}\ar@{.}[uuuullll]& & *{a_2} \ar@{-}[uurr] \ar@{.}[dddrrr]& & && *{}\\
    &&&*{}& & & & *{} \ar@{.}[uuuullll] & & *{a_1} \ar@{.}[ddrr] \ar@{-}[ur]& & &{\!\!\!\! \ell_2} &*{}&&&&&&*{}\\ 
    &&& & & & & & & {\;\;\;\; \ell_1} &&&& \\
    *{}\ar@{--}[r]&*{} \ar@{-}[rrrr] &&& & *{} \ar@{.}[uuuullll]\ar@{.}[uuuuuuuuurrrrrrrrr]\ar@{-}[rr]& *{} & *{} \ar@{-}[uurr]\ar@{-}[rrrr]& & & & *{} \ar@{.}[uuuuuuurrrrrrr]\ar@{-}[rrrrrrrrrr]&&\ar@{.}[uuuuuurrrrrr]&*{}&*{}&&*{}\ar@{.}[uuuurrrr]&&&&*{}\ar@{--}[r]&*{}\\
           }
}
\end{equation}
Note that the regions $H^-(\Sigma a)$ for $a$ in $P$ share the half
line $s$ as a common edge, while each of the line segments $\ell_i$
is an edge of a region $H^+(\Sigma a)$ with $a$ in $P$.

We are aiming to show that $n$ is also a left-fountain, that is, there
are infinitely many arcs in $\arcs$ of the form $(m,n)$.  This is the
same as showing that there are infinitely many indecomposable objects
of $\subcat$ which are on the half line $s$.

Let $x$ be an indecomposable object in the region bounded by the
diagonal half lines $s$ and $s^{\prime}$ and let $b \rightarrow x$ be
an $\subcat$-precover.  We can assume that the morphism $b \rightarrow
x$ is non-zero on each direct summand of $b$; in particular, each
direct summand of $b$ belongs to $H(\Sigma^{-1}x)$.

It is easy to see that $x$ is in $H^-(\Sigma a)$ for infinitely many
$a$ in $P$, so there are infinitely many $a$ in $P$ with a non-zero
morphism $a \rightarrow x$.  Each of these morphisms factors through
$b \rightarrow x$, so there is an indecomposable direct summand $c$ of
$b$ such that infinitely many $a$ in $P$ have non-zero morphisms to
$c$.  Hence $c$ is certainly in $H(\Sigma a)$ for some $a$ in $P$.
Moreover, since $c$ is in $\subcat$, Remark \ref{rmk:forbidden} says
that $c$ must be outside the forbidden region $H(a)$ for each $a$ in
$P$.  The following sketch shows the regions $H(\Sigma a)$ (ordinary
lines) and $H(a)$ (wavy lines) for an indecomposable object $a$.
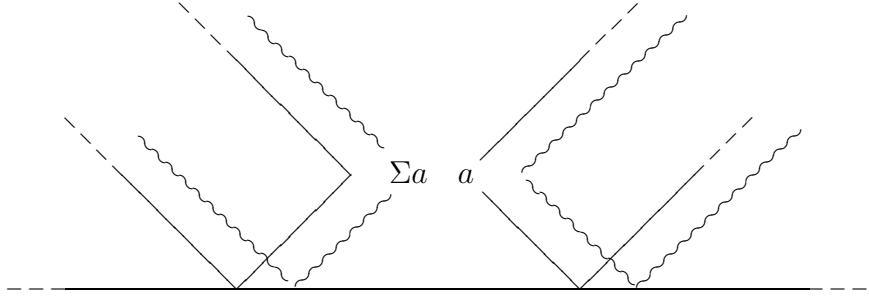
\[
  \xymatrix @-2.0pc @! {
    &&&*{} &&&&&&&& *{}&& \\
    &&&& *{} \ar@{--}[ul] & & & & & & *{} \ar@{--}[ur] \\
    &*{}&& & & & & & & & & && *{}&\\
    &&*{}\ar@{--}[ul]& & & & *{} \ar@{-}[ddll] \ar@{-}[uull] & \Sigma a\ar@{~}[uuulll]\ar@{~}[ddll]& a \ar@{-}[ddrr] \ar@{-}[uurr]& *{}\ar@{~}[uuurrr]\ar@{~}[ddrr]& &&*{}\ar@{--}[ur]&\\ 
    && \\
    *{}\ar@{--}[r]&*{} \ar@{-}[rrr] && & *{} \ar@{-}[uull]\ar@{-}[rrrrrr]& *{}\ar@{~}[uuulll]& & & & & *{} \ar@{-}[uurr]\ar@{-}[rrrr]&*{}\ar@{~}[uuurrr]&&&*{}\ar@{--}[r]&*{}\\
           }
\]
Combining this with the previous sketch shows that there are only
three possible places for $c$: It is either on one of the line
segments $\ell_i$, or on the half line $r$, or on the half line $s$.

Now, there are infinitely many $a$ in $P$ with non-zero morphisms to
$c$; that is, in\-fi\-ni\-te\-ly many $a$ in $P$ which are in
$H(\Sigma^{-1}c)$.  And it is easy to see from the sketch
\eqref{equ:fig1} that this does not happen if $c$ is on an $\ell_i$ or
on $r$, so $c$ must be on $s$.  Moreover, $c$ is a direct summand of
$b$, so $c$ is in $H(\Sigma^{-1}x)$.  Combining the sketch
\eqref{equ:fig1} with $H(\Sigma^{-1}x)$, indicated in wavy lines,
gives the following.
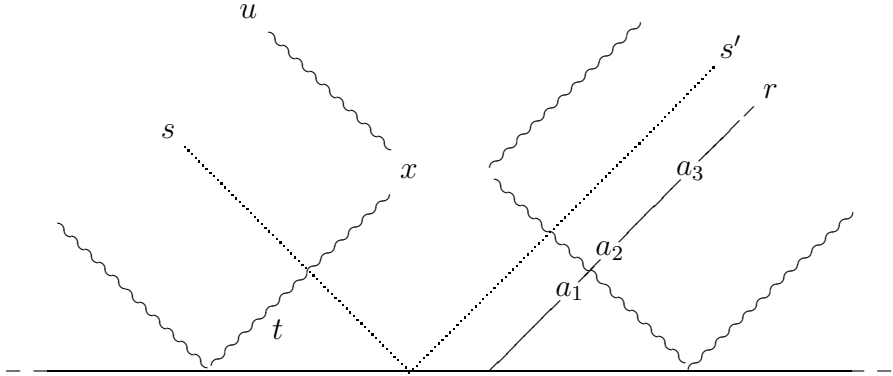
\[
\vcenter{
  \xymatrix @-2.05pc @! {
    &&&&&&{u}&&&&&*{}&&&&&&&&&& \\
    &&&&&&&&&&&&&&&&&&{s^{\prime}}&& \\
    &&&&&&&&&*{}&&&&&&&&&&{r}&&&&*{}\\
    &&&&{s}&&&&*{}&&&&&&&&&&*{}\ar@{--}[ur]&&&&\\
    &&&&&&&&&& {x} \ar@{~}[dddddlllll]\ar@{~}[uuuullll]&&*{}\ar@{~}[dddddrrrrr]\ar@{~}[uuuurrrr]&&*{}& *{} & *{} & *{a_3} \ar@{-}[ur] &&&&*{}& \\
    &&&&&&&& &  & & &&& & & *{} &&&&*{}&*{}&*{}\\
    &&&&&&&*{}&& & & & & *{}& & *{a_2} \ar@{-}[uurr] & & && *{}&&&\\
    &&&&&&&&*{}& & & & *{}  & & *{a_1} \ar@{-}[ur]& & & &*{}&&&&*{}\\ 
    &&&& & & & {\!\!\!\! t} &&&&& & & &&&& \\
    *{}\ar@{--}[r]&*{}\ar@{-}[rrrr]&&&&*{} \ar@{~}[uuuullll]\ar@{-}[rrrrr] &&&& & *{} \ar@{.}[uuuuuullllll]\ar@{.}[uuuuuuuurrrrrrrr]\ar@{-}[rr]& *{} & *{} \ar@{-}[uurr]\ar@{-}[rrrr]& & & & *{} \ar@{-}[rrrrr]&*{}\ar@{~}[uuuurrrr]&&*{}&&*{}\ar@{--}[r]&*{}\\
           }
}
\]
This shows that the indecomposable object $c$ must be on $s$ and above
the line segment $t$.

For each $x$ we get a $c$ in $\subcat$ which is on $s$ and above the
line segment $t$ corresponding to $x$.  By moving $x$ out along the
half line $u$, we can clearly force infinitely many distinct $c$'s.
It follows that, as desired, there are infinitely many indecomposable
objects of $\subcat$ which are on $s$.

Suppose next that the right-fountain of $\arcs$ is also a
left-fountain.  Let $x$ be an indecomposable object of $\sD$; we will
construct a $\subcat$-precover of $x$.  Consider again the sketch
\eqref{equ:fig1}.  The arcs going right, respectively left from the
fountain of $\arcs$ correspond to in\-de\-com\-po\-sa\-ble objects of
$\subcat$ on the half lines $r$, respectively $s$, so each of $r$ and
$s$ contains infinitely many indecomposable objects of $\subcat$.  The
other arcs in $\arcs$ correspond to indecomposable objects of
$\subcat$ away from $r$ and $s$.

Consider the set of indecomposable objects $a$ of $\subcat$ which have
non-zero morphisms to $x$, and divide it into disjoint subsets $R$,
$S$, and $T$ according to whether $a$ is on $r$, $s$, or neither.  We
will construct a morphism $a_r \rightarrow x$ with $a_r \in \subcat$
such that each $a \rightarrow x$ with $a \in R$ factors as $a
\rightarrow a_r \rightarrow x$.  We will also construct morphisms $a_s
\rightarrow x$ and $a_t \rightarrow x$ with the analogous properties
with respect to $S$ and $T$; then an $\subcat$-precover can be
obtained as $a_r \oplus a_s \oplus a_t \rightarrow x$.

If a set $R$, $S$, or $T$ is finite, then the construction of the
corresponding morphism $a_r \rightarrow x$, $a_s \rightarrow x$, or
$a_t \rightarrow x$ is trivial.

The set $T$ is always finite: Suppose that $a$ is in $T$ and let $\fa$
be the arc corresponding to $a$. There is a non-zero morphism $a
\rightarrow x$ so Lemma \ref{lem:Hom} gives that $\fa$ intersects the
arc $\fx = (i,j)$ corresponding to $\Sigma^{-1}x$.  Hence $\fa$ ends
in an integer $m$ with $i < m < j$.  Since $a$ is in $T$, it is in
$\subcat$ but not on one of the half lines $r$ and $s$; this means
that $\fa$ is an arc which is in $\arcs$ but does not end in the
fountain of $\arcs$.  In particular, $m$ is not the fountain of
$\arcs$.  We conclude that each of the finitely many possible values
of $m$ is an integer where only finitely many arcs of $\arcs$ end, and
it follows that there are only finitely many arcs $\fa$ as described.
That is, $T$ has finitely many elements.

We are left to deal with the cases of $R$ and $S$ being infinite.

Suppose that $R$ is infinite.  So there are infinitely many
indecomposable objects $a$ of $\subcat$ on the half line $r$ with
non-zero morphisms to $x$, that is, with $x$ in $H(\Sigma a)$.  By
inspecting the sketch \eqref{equ:fig1} it can be seen that $x$ is in
the region bounded by the half lines $s$ and $s^{\prime}$.  However,
there are infinitely many indecomposable objects of $\subcat$ on the
half line $s$, and so we can chose one, $a_r$, which has $x \in
H^+(\Sigma a_r)$ as indicated here.
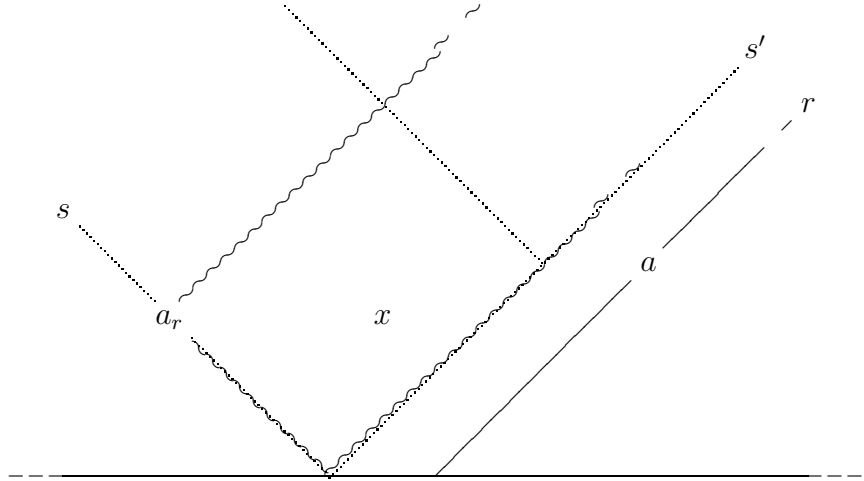
\[
  \xymatrix @-1.8pc @! {
    &&&&&&&*{}&&&&&&&&& \\
    &&&&&&&&*{}\ar@{~~}[ur]&&&&&&{s^{\prime}}& \\
    &&&&&*{}&&&&&&&&*{}&&{r}&\\
    &&&&*{}&&&&&&&&&&*{}\ar@{--}[ur]&&\\
    &{s}&&&&*{} &&&&&*{}&*{}\ar@{~~}[ur]& *{} & &&& \\
    &&&&&& & & & & *{} \ar@{.}[uuuuulllll]& & {a} \ar@{-}[uurr] &&&&*{}\\
    &&& {a_r} \ar@{.}[uull] \ar@{~}[uuuuurrrrr]&& & & {x}& & & & & & && *{}\\
    &&&&*{}& & & & & & & & & &*{}&&\\ 
    &&&& & & & & & & &&&& \\
    *{}\ar@{--}[r]&*{} \ar@{-}[rrrrr] &&&& & *{} \ar@{.}[uuulll]\ar@{.}[uuuuuuuurrrrrrrr]\ar@{~}[uuulll]\ar@{~}[uuuuurrrrr]\ar@{-}[rr]& *{} & *{} \ar@{-}[uuuurrrr]\ar@{-}[rrrr]& & & & *{} \ar@{-}[rrr]&&*{}&*{}\ar@{--}[r]&*{}\\
           }
\]
Pick a non-zero morphism $a_r \rightarrow x$.  If $a$ is in $R$ then
it has a non-zero morphism $a \rightarrow x$, and then we have $x \in
H^-(\Sigma a)$ as in the sketch.  But it is clear that $a_r \in
H^-(\Sigma a)$ and so Lemma \ref{lem:xyz2} says that a morphism $a
\rightarrow x$ factors like $a \rightarrow a_r \rightarrow x$ as
desired.

Suppose that $S$ is infinite.  So there are infinitely many
indecomposable objects $a$ of $\subcat$ on the half line $s$ with
non-zero morphisms to $x$, that is, with $x$ in $H(\Sigma a)$.  The
following sketch shows some of the indecomposable objects $a$ on $s$
and, in dotted lines, their regions $H(\Sigma a)$.  Since $x$ is in
infinitely many of these regions, it can be seen that it is again in
the region bounded by the half lines $s$ and $s^{\prime}$.
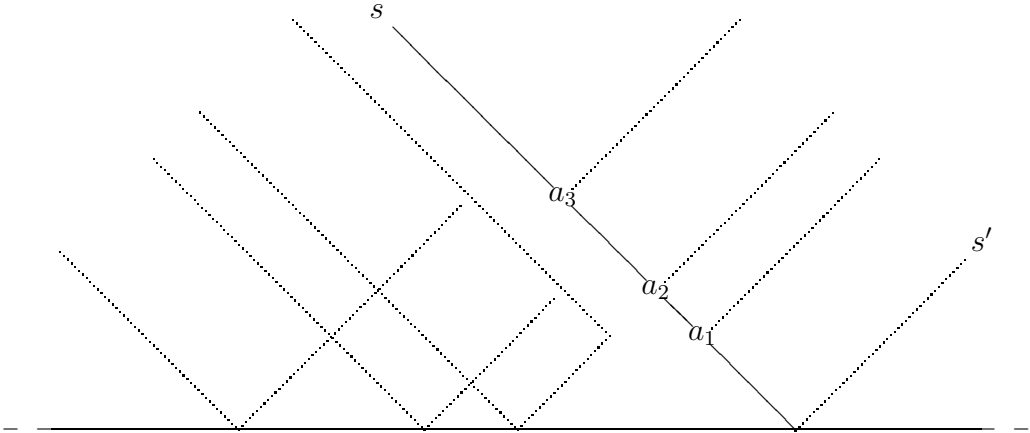
\[
%  \def\objectstyle{\scriptscriptstyle}
%  \xymatrix @-2.05pc @! {
  \xymatrix @-1.85pc @! {
    &&&&&&&&{s}&&&&&&&&& \\
    &&&&&&&&&&&&&&&&&&& \\
    &&&&&&&&&&&&&&&&&&& \\
    &&&&&&&&&&&&&&&&&&& \\
    &&&&&&&&&&&&*{a_3}\ar@{-}[uuuullll]\ar@{.}[uuuurrrr]&&&&&&& \\
    &&&&&&&&&&&&&&&&&&&&&{s^{\prime}} \\
    &&&&&&&&&&&&&&*{a_2}\ar@{-}[uull]\ar@{.}[uuuurrrr]&&&&& \\
    &&&&&&&&&&&&&*{}\ar@{.}[uuuuuuulllllll]&&*{a_1}\ar@{-}[ul]\ar@{.}[uuuurrrr]&&&& \\
    &&&&&&&&&&&&&&&&&&& \\
    *{}\ar@{--}[r]&*{}\ar@{-}[rrrrrrrrrrrrrrrrrrrr]&&&&\ar@{.}[uuuullll]\ar@{.}[uuuuurrrrr]&&&&\ar@{.}[uuuuuullllll]\ar@{.}[uuurrr]&&\ar@{.}[uuuuuuulllllll]\ar@{.}[uurr]&&&&&&*{}\ar@{-}[uull]\ar@{.}[uuuurrrr]&&&&*{}\ar@{--}[r]&*{} \\
           }
\]
Let $a_s$ be the indecomposable object in $S$ which is closest to the
end of $s$ and pick a non-zero morphism $a_s \rightarrow x$.  It is
clear that we have $x \in H^+(\Sigma a_s)$.  If $a$ is in $S$ then it has
a non-zero morphism $a \rightarrow x$, and again $x \in H^+(\Sigma
a)$.  The following sketch shows the whole situation.
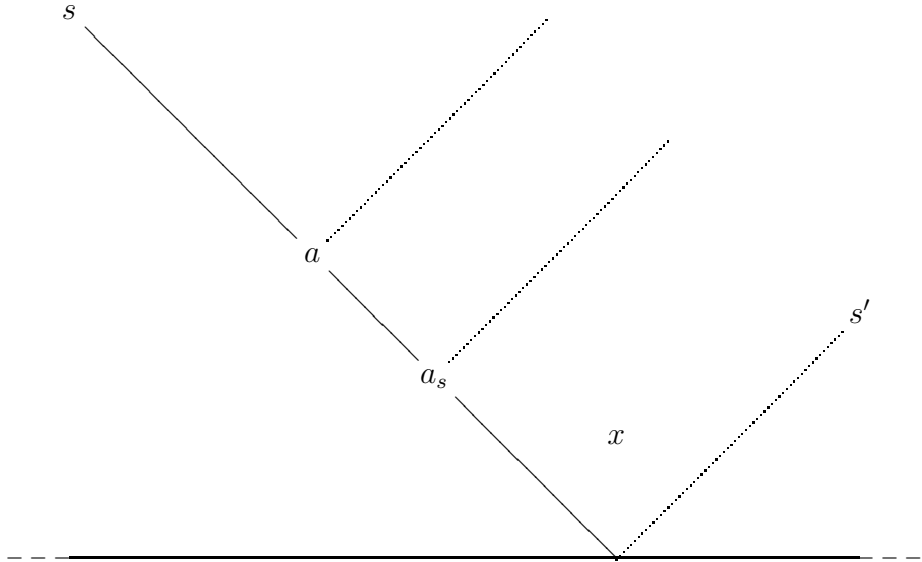
\[
%  \def\objectstyle{\scriptscriptstyle}
%  \xymatrix @-1.7pc @! {
  \xymatrix @-1.55pc @! {
    &{s}&&&&&&&&& \\
    &&&&&&&&&&&& \\
    &&&&&&&&&&&& \\
    &&&&&&&&&&&& \\
    &&&&&{a}\ar@{-}[uuuullll]\ar@{.}[uuuurrrr]&&&&&&& \\
    &&&&&&&&&&&&&&{s^{\prime}} \\
    &&&&&&&{a_s}\ar@{-}[uull]\ar@{.}[uuuurrrr]&&&&& \\
    &&&&&&&&&& {x} && \\
    &&&&&&&&&&&& \\
    *{}\ar@{--}[r]&*{}\ar@{-}[rrrrrrrrrrrrr]&&&&&&&&&*{}\ar@{-}[uuulll]\ar@{.}[uuuurrrr]&&&&*{}\ar@{--}[r]&*{} \\
           }
\]
But it is clear that $a_s \in H^+(\Sigma a)$ and so Lemma
\ref{lem:xyz1}(ii) says that each morphism $a \rightarrow x$ factors
like $a \rightarrow a_s \rightarrow x$ as desired.
\end{proof}

\begin{Example}
\label{exa:non-functorially_finite}
Theorem \ref{thm:cluster_tilting} shows that there are maximal
$1$-orthogonal subcategories of $\sD$ which are not cluster tilting;
that is, they are not functorially finite.

A concrete example comes from the maximal set of non-intersecting arcs
in the sketch \eqref{equ:right_and_left_fountain} in the Introduction,
which corresponds to the maximal $1$-orthogonal subcategory $\subcat$
with the indecomposable objects marked by bullets.
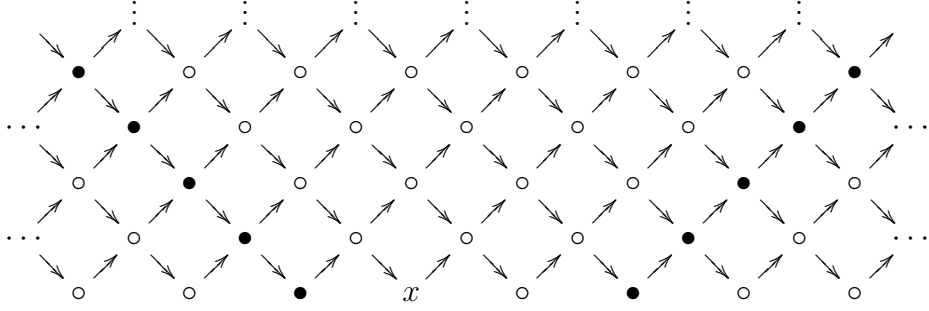
\[
  \xymatrix @-2.1pc @! {
    \rule{2ex}{0ex} \ar[dr] & & \vdots \ar[dr]  & & \vdots \ar[dr]& & \vdots \ar[dr]& & \vdots \ar[dr] & & \vdots \ar[dr] & & \vdots \ar[dr] & & \vdots \ar[dr] & & \rule{2ex}{0ex} \\
    & \bullet \ar[ur] \ar[dr] & & \circ \ar[ur] \ar[dr] & &\circ \ar[ur] \ar[dr]  & & \circ \ar[ur] \ar[dr]& & \circ \ar[ur] \ar[dr] & & \circ \ar[ur] \ar[dr] & & \circ \ar[ur] \ar[dr] & & \bullet \ar[ur] \ar[dr] & \\
    \cdots \ar[ur]\ar[dr]& & \bullet \ar[ur] \ar[dr] & & \circ \ar[ur] \ar[dr] &  & \circ \ar[ur] \ar[dr] && \circ \ar[ur] \ar[dr] & & \circ \ar[ur] \ar[dr] & & \circ \ar[ur] \ar[dr] & & \bullet \ar[ur] \ar[dr] & & \cdots \\
    & \circ \ar[ur] \ar[dr] & & \bullet \ar[ur] \ar[dr] & &\circ \ar[ur] \ar[dr] & &\circ \ar[ur] \ar[dr] & & \circ \ar[ur] \ar[dr] & & \circ \ar[ur] \ar[dr] & & \bullet \ar[ur] \ar[dr] & & \circ \ar[ur] \ar[dr] & \\
    \cdots \ar[ur]\ar[dr]& & \circ \ar[ur] \ar[dr] &  & \bullet \ar[ur] \ar[dr] && \circ \ar[ur] \ar[dr]& & \circ \ar[ur] \ar[dr] & & \circ \ar[ur] \ar[dr] & & \bullet \ar[ur] \ar[dr] & & \circ \ar[ur] \ar[dr] & & \cdots\\
    & \circ \ar[ur] & & \circ \ar[ur] & &\bullet \ar[ur] & & x \ar[ur] & & \circ \ar[ur] & & \bullet \ar[ur] & & \circ \ar[ur] & & \circ \ar[ur] & \\
               }
\]
In fact, it is not hard to adapt the arguments in the proof of Theorem
\ref{thm:cluster_tilting} to show that there is no $\subcat$-precover
of the indecomposable object $x$.
\end{Example}

\section{The cluster structure of $\sD$}
\label{sec:cluster}

This section proves Theorem C from the Introduction; see Theorem
\ref{thm:cluster}. 

\begin{Definition}
For each cluster tilting subcategory $\subcat$ of $\sD$ we can
consider the set $\cluster$ of indecomposable objects of $\subcat$
whence $\subcat = \add \cluster$.  We will refer to the sets
$\cluster$ as {\em clusters}.
\end{Definition}

The clusters are said to form a {\em cluster structure} if the
following conditions are satisfied; cf.\ \cite{BIRS}.

\begin{enumerate}

  \item If $\cluster$ is a cluster, then each of its indecomposable
  objects $a$ can be replaced with a unique other indecomposable
  object $a^*$ of $\sD$ such that a new cluster $\cluster^*$
  results.

\smallskip

  \item There are distinguished triangles $a^* \rightarrow b \rightarrow
  a$ and $a \rightarrow b^{\prime} \rightarrow a^*$ in $\sD$ where the
  left-hand morphisms are $\add(\cluster \! \setminus \! \{ a
  \})$-envelopes and the right-hand morphisms are $\add(\cluster \!
  \setminus \! \{ a \})$-covers.

\smallskip

  \item If $\cluster$ is a cluster, then the AR quiver of $\add
  \cluster$ has no loops or $2$-cycles.

\smallskip

  \item Passing from the AR quiver of $\add \cluster$ to the AR quiver
  of $\add \cluster^*$ is given by Fomin-Zelevinsky mutation at $a$
  in the sense of \cite[sec.\ 8]{FZ2}.

\end{enumerate}

\begin{Theorem}
\label{thm:cluster}
The clusters form a cluster structure on $\sD$.
\end{Theorem}

\begin{proof}
Remark \ref{rmk:basics} says that $\sD$ is a $2$-Calabi-Yau category
and it follows from Theorem \ref{thm:cluster_tilting} that there exist
cluster tilting subcategories of $\sD$.  Hence by \cite[thm.\
I.1.6]{BIRS} it is enough to show that for each cluster $\cluster$,
there are no loops or $2$-cycles in the AR quiver of the cluster
tilting subcategory $\subcat = \add \cluster$.

It is clear that there are no loops: If $a$ is in $\cluster$, then
$\subcat(a,a) = \sD(a,a) = k$ by Corollary \ref{cor:to_and_from}, so
each non-zero morphism $a \rightarrow a$ is an isomorphism and so not
irreducible.

To show that there are no $2$-cycles in the AR quiver of $\subcat =
\add \cluster$, we will show the stronger claim that given $a$ and
$b$ in $\cluster$ with $\subcat(a,b) \neq 0$, it follows that
$\subcat(b,a) = 0$.

The following sketch shows the regions $H(\Sigma a)$ (straight lines)
and $H(a)$ (wavy lines).
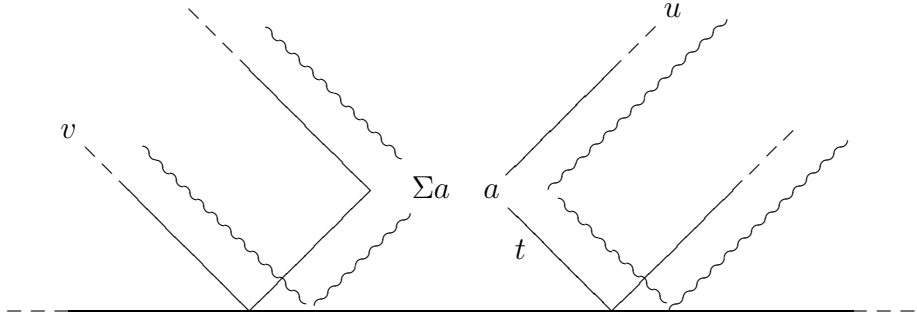
\[
  \xymatrix @-1.9pc @! {
    &&&*{} &&&&&&&& {u}&& \\
    &&&& *{} \ar@{--}[ul] & & & & & & *{} \ar@{--}[ur] \\
    &{v}&& & & & & & & & & && *{}&\\
    &&*{}\ar@{--}[ul]& & & & *{}\ar@{-}[ddll] \ar@{-}[uull] & {\Sigma a}\ar@{~}[uuulll]\ar@{~}[ddll]& {a} \ar@{-}[ddrr] \ar@{-}[uurr]& *{}\ar@{~}[uuurrr]\ar@{~}[ddrr]& &&*{}\ar@{--}[ur]&\\ 
    &&&*{} &&&&&&{\!\!\!\!\!\!\!\!\!\!\!\! t} \\
    *{}\ar@{--}[r]&*{} \ar@{-}[rrr] && & *{} \ar@{-}[uull]\ar@{-}[rrrrrr]& *{}\ar@{~}[uuulll]& & & & & *{} \ar@{-}[uurr]\ar@{-}[rrrr]&*{}\ar@{~}[uuurrr]&&&*{}\ar@{--}[r]&*{}\\
           }
\]
Since $\sD(a,b) = \subcat(a,b) \neq 0$, the indecomposable object $b$
is in the region $H(\Sigma a)$ by Corollary
\ref{cor:to_and_from}.  Since $a$ and $b$ are both in the cluster
tilting subcategory $\subcat$, the object $b$ is outside the region
$H(a)$ by Remark \ref{rmk:forbidden}.

It follows that in the sketch, $b$ must be either on the line segment
$t$ or on one of the half lines $u$ and $v$, and in any of these cases
it is easy to verify that $a$ is outside $H(\Sigma b)$, that is,
$\subcat(b,a) = \sD(b,a) = 0$.
\end{proof}

\begin{Example}
\label{exa:mutation}
Passing from the cluster $\cluster$ to $\cluster^*$ is referred to as
cluster mutation at $a$.  It corresponds to an obvious combinatorial
mutation of maximal sets of arcs.

For instance, recall the leapfrog configuration.
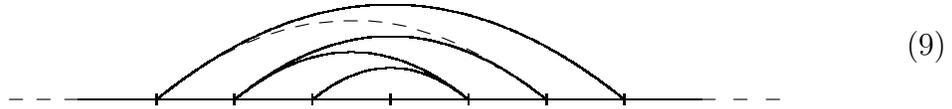
\begin{equation}
\label{equ:fig2}
\vcenter{
  \xymatrix @-4.0pc @! {
       \rule{0ex}{9ex} \ar@{--}[r]
     & *{}\ar@{-}[r]
     & *{\rule{0.1ex}{0.8ex}} \ar@{-}[r] \ar@/^2.5pc/@{--}[rrrrr]\ar@/^3.0pc/@{-}[rrrrrr]
     & *{\rule{0.1ex}{0.8ex}} \ar@{-}[r] \ar@/^1.5pc/@{-}[rrr]\ar@/^2.0pc/@{-}[rrrr]
     & *{\rule{0.1ex}{0.8ex}} \ar@{-}[r] \ar@/^1.0pc/@{-}[rr]
     & *{\rule{0.1ex}{0.8ex}} \ar@{-}[r] 
     & *{\rule{0.1ex}{0.8ex}} \ar@{-}[r]
     & *{\rule{0.1ex}{0.8ex}} \ar@{-}[r]
     & *{\rule{0.1ex}{0.8ex}} \ar@{-}[r]
     & *{}\ar@{--}[r]
     & *{}
                    }
}
\end{equation}
Under the bijection of Remark \ref{rmk:bijection}, this corresponds to
the following cluster $\cluster$ in $\sD$; the broken arc corresponds
to the object $a$.
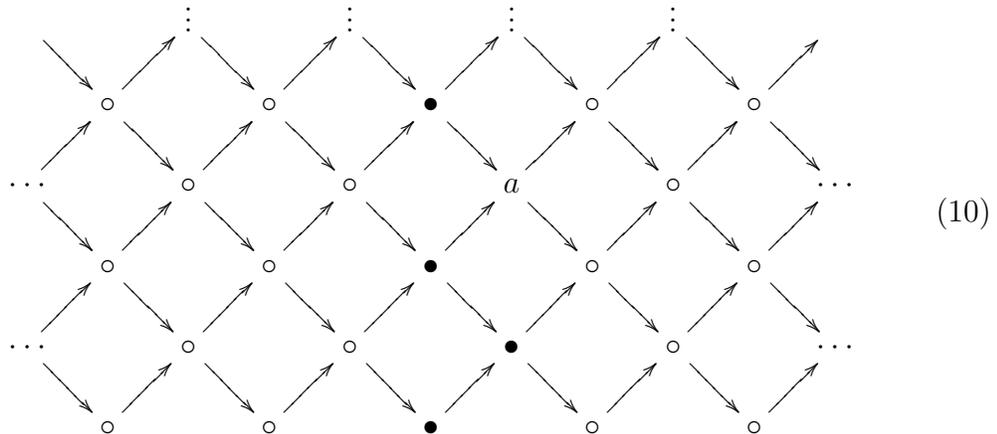
\begin{equation}
\label{equ:A}
\vcenter{
  \xymatrix @-1.3pc @! {
    \rule{2ex}{0ex} \ar[dr] & & \vdots \ar[dr] & & \vdots \ar[dr] & & \vdots \ar[dr] & & \vdots \ar[dr] & & \rule{2ex}{0ex} \\
    & \circ \ar[ur] \ar[dr] & & \circ \ar[ur] \ar[dr] & & \bullet \ar[ur] \ar[dr] & & \circ \ar[ur] \ar[dr] & & \circ \ar[ur] \ar[dr] &  \\
    \cdots \ar[ur]\ar[dr]& & \circ \ar[ur] \ar[dr] & & \circ \ar[ur] \ar[dr] & & a \ar[ur] \ar[dr] & & \circ \ar[ur] \ar[dr] & & \cdots \\
    & \circ \ar[ur] \ar[dr] & & \circ \ar[ur] \ar[dr] & & \bullet \ar[ur] \ar[dr] & & \circ \ar[ur] \ar[dr] & & \circ \ar[ur] \ar[dr] & \\
    \cdots \ar[ur]\ar[dr]& & \circ \ar[ur] \ar[dr] & & \circ \ar[ur] \ar[dr] & & \bullet \ar[ur] \ar[dr] & & \circ \ar[ur] \ar[dr] & & \cdots \\
    & \circ \ar[ur] & & \circ \ar[ur] & & \bullet \ar[ur] & & \circ \ar[ur] & & \circ \ar[ur] & \\
               }
}
\end{equation}
Removing the broken arc from \eqref{equ:fig2} creates a `quadrangle',
and there is clearly a unique other arc which bisects it to form a new
maximal set of non-intersecting arcs.
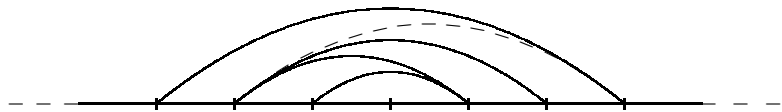
\[
\vcenter{
  \xymatrix @-4.0pc @! {
       \rule{0ex}{9ex} \ar@{--}[r]
     & *{}\ar@{-}[r]
     & *{\rule{0.1ex}{0.8ex}} \ar@{-}[r] \ar@/^3.0pc/@{-}[rrrrrr]
     & *{\rule{0.1ex}{0.8ex}} \ar@{-}[r] \ar@/^2.5pc/@{--}[rrrrr]\ar@/^1.5pc/@{-}[rrr]\ar@/^2.0pc/@{-}[rrrr]
     & *{\rule{0.1ex}{0.8ex}} \ar@{-}[r] \ar@/^1.0pc/@{-}[rr]
     & *{\rule{0.1ex}{0.8ex}} \ar@{-}[r] 
     & *{\rule{0.1ex}{0.8ex}} \ar@{-}[r]
     & *{\rule{0.1ex}{0.8ex}} \ar@{-}[r]
     & *{\rule{0.1ex}{0.8ex}} \ar@{-}[r]
     & *{}\ar@{--}[r]
     & *{}
                    }
}
\]
Under the bijection of Remark \ref{rmk:bijection}, this corresponds to
the cluster $\cluster^*$; the broken arc corresponds to the object
$a^*$.
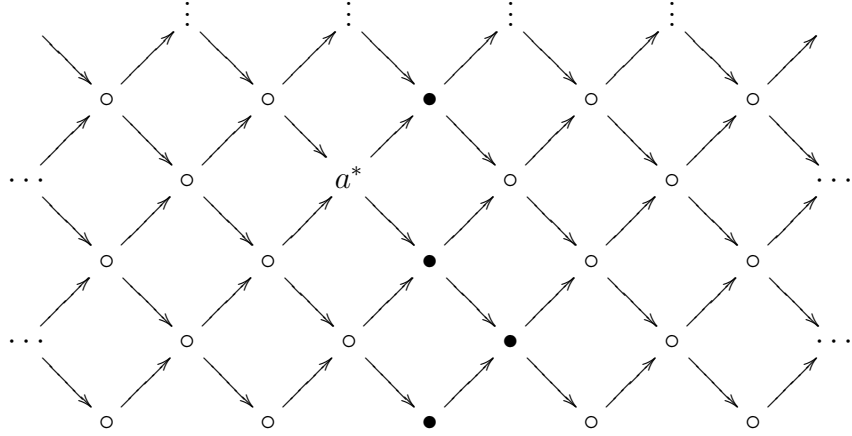
\[
  \xymatrix @-1.3pc @! {
    \rule{2ex}{0ex} \ar[dr] & & \vdots \ar[dr] & & \vdots \ar[dr] & & \vdots \ar[dr] & & \vdots \ar[dr] & & \rule{2ex}{0ex} \\
    & \circ \ar[ur] \ar[dr] & & \circ \ar[ur] \ar[dr] & & \bullet \ar[ur] \ar[dr] & & \circ \ar[ur] \ar[dr] & & \circ \ar[ur] \ar[dr] &  \\
    \cdots \ar[ur]\ar[dr]& & \circ \ar[ur] \ar[dr] & & a^* \ar[ur] \ar[dr] & & \circ \ar[ur] \ar[dr] & & \circ \ar[ur] \ar[dr] & & \cdots \\
    & \circ \ar[ur] \ar[dr] & & \circ \ar[ur] \ar[dr] & & \bullet \ar[ur] \ar[dr] & & \circ \ar[ur] \ar[dr] & & \circ \ar[ur] \ar[dr] & \\
    \cdots \ar[ur]\ar[dr]& & \circ \ar[ur] \ar[dr] & & \circ \ar[ur] \ar[dr] & & \bullet \ar[ur] \ar[dr] & & \circ \ar[ur] \ar[dr] & & \cdots \\
    & \circ \ar[ur] & & \circ \ar[ur] & & \bullet \ar[ur] & & \circ \ar[ur] & & \circ \ar[ur] & \\
               }
\]
\end{Example}

\begin{Example}
\label{exa:Gamma}
If $\subcat$ is a cluster tilting subcateory of $\sD$, then
$\sD/\subcat$ is an abelian category by \cite[sec.\ 2]{KellerReiten2}
and \cite[thm.\ 3.3]{KoenigZhu}, and we have $\sD/\subcat \simeq
\mod\,\subcat$ by \cite[sec.\ 2]{KellerReiten2} and \cite[cor.\
4.4]{KoenigZhu}.

In the case of the $\subcat$ given by the sketch \eqref{equ:A} above,
it is easy to check that $\cA$ is the path category of its AR quiver
$\Gamma$.
\[
  \Gamma = 
  \xymatrix @! {
    \bullet \ar[r] & \bullet & \bullet \ar[l]\ar[r] & a & \bullet \ar[l]\ar[r] & \cdots
               }
\]
So $\mod\,\subcat$ is equivalent to $\rep \Gamma$, the category of
finitely presented representations of $\Gamma$.  Such hereditary
categories were studied in \cite{RVdB}.

Likewise, $\cA^*$ is the path category of its AR quiver $\Gamma^*$.
\[
  \Gamma^* = 
  \xymatrix @! {
    \bullet \ar[r] & \bullet & \bullet \ar[l] & a^* \ar[l]\ar[r]& \bullet\ar[r] & \cdots
               }
\]
So $\mod\,\subcat^*$ is equivalent to $\rep \Gamma^*$, and cluster
mutation at $a$ has changed $\rep \Gamma$ to $\rep \Gamma^*$. 
\end{Example}

\section{Questions}
\label{sec:questions}

Let us end the paper by posing some questions which seem natural in
the light of the results presented here.

(1) The category $\sD$ behaves like a cluster category of type
$A_{\infty}$.  Is it possible to define a cluster algebra of type
$A_{\infty}$?

(2)  Section \ref{sec:cluster} gives the means to do cluster tilting
of abelian categories of the form $\sD/\subcat$ where $\subcat$ is a
cluster tilting subcategory.  Which abelian categories have this form?
In particular, which hereditary abelian categories do?

(3)  Can $\sD$ be viewed as a covering category for the tubular
$2$-Calabi-Yau categories studied in \cite[sec.\ 2]{BMV}?

(4)  The AR quiver of $\sD$ is ${\mathbb Z}A_{\infty}$.  Is there a
similar category with AR quiver ${\mathbb Z}\Delta$ when $\Delta$ is
another infinite Dynkin quiver than $A_{\infty}$?

(5)  Is it possible to define `higher cluster categories of type
$A_{\infty}$'?  See \cite{Keller}, \cite{Thomas}, and \cite{BinZhu} for
the type $A_n$ case.

\medskip
\noindent
{\bf Acknowledgement. }
We thank Collin Bleak and Robert Marsh for interesting discussions on
the material presented here.  In particular, Collin Bleak showed us
that if a triangulation of the $\infty$-gon has a right-fountain then
it also has a left-fountain, and Robert Marsh explained to us how to
use the tools of \cite{BIRS} to prove Theorem C.  We thank Katsuhiko
Kuribayashi for pointing out that $S^2$ is formal over any field, and
Bernard Leclerc, Robert Marsh, and Andrei Zelevinsky for remarks
on a preliminary version.


\begin{thebibliography}{19}
  
% \bibitem{Amiot} C.\ Amiot, {\it On the structure of triangulated
%     categories with finitely many in\-de\-com\-po\-sa\-bles}, Bull.\
%   Soc.\ Math.\ France, to appear.  {\tt math.CT/0612141}.
  
% \bibitem{Asashiba} H.\ Asashiba, {\it On a lift of an individual
% stable equivalence to a standard derived e\-qui\-va\-len\-ce for
% representation-finite self-injective algebras}, Algebr.\ 
% Represent.\ Theory\ {\bf 6} (2003), 427--447.

%\bibitem{ABST}  I.\ Assem, T.\ Br\"ustle, R.\ Schiffler, G.\ Todorov,
%$m$-cluster categories and $m$-replicated algebras,
%preprint (2006). {\tt math.RT/0608727}.

% \bibitem{ASS}  I.\ Assem, D.\ Simson, A.\ Skowro\'{n}ski, ``Elements of
% the Representation Theory of Associative Algebras'', Vol.\ 1, London
% Math.\ Soc.\ Stud.\ Texts, Vol.\ 65, Cambridge University Press,
% Cambridge, 2006.

%\bibitem{AusArtI}  M.\ Auslander, {\it Representation theory of Artin
%algebras I}, Comm.\ Algebra {\bf 1} (1974), 177--268.

%\bibitem{AusArtIII}  \bysame I.\ Reiten, {\it Representation
%theory of Artin algebras III}, Comm.\ Algebra {\bf 3} (1975),
%239--294. 

%\bibitem{ARS}  \bysame, I.\ Reiten, S.\ O.\ Smal\o,
%``Representation theory of Artin algebras'', Cambridge Stud.\ Adv.\
%Math., Vol.\ 36, Cambridge University Press, Cambridge, 1997, first
%paperback edition with corrections.

%\bibitem{BaurMarsh1} K.\ Baur, R.\ Marsh, {\it A geometric description
%of $m$-cluster categories}, to appear in Trans.\ Amer.\ Math.\ 
%Soc.

%\bibitem{BaurMarsh2}  K.\ Baur, R.\ Marsh, {\it A geometric
%description of the $m$-cluster categories of type $D_n$}, to appear in
%Int.\ Math.\ Res.\ Not.

% \bibitem{BS}  J.\ Bia{\l}kowski, A.\ Skowro\'{n}ski, 
% {\it Calabi-Yau stable module categories of finite type}, 
% Colloq.\ Math.\ {\bf 109} (2007), 257--269. 

\bibitem{BIRS}  A.\ B.\ Buan, O.\ Iyama, I.\ Reiten, and J.\ Scott,
{\it Cluster structures for $2$-Calabi-Yau categories and unipotent
groups}, to appear in Compositio Math.  {\tt math.RT/0701557.}


\bibitem{BMRRT}  A.\ B.\ Buan, R.\ J.\ Marsh, M.\ Reineke, I.\ Reiten,
and G.\ Todorov, {\it Tilting theory and cluster combinatorics}, Adv.\
Math.\ {\bf 204} (2006), 572--618. 

% \bibitem{BMR} A.\ B.\ Buan, R.\ J.\ Marsh, I. Reiten, {\it
% Cluster-tilted algebras}, Trans.\ Amer.\ Math.\ Soc.\ {\bf 359}
% (2007), 323--332.

\bibitem{BMV}  A.\ B.\ Buan, R.\ J.\ Marsh, and D.\ F.\ Vatne, Cluster
structures from $2$-Calabi-Yau categories with loops, preprint
(2008).  {\tt math.RT/0810.3132}.

\bibitem{CCS}  P.\ Caldero, F.\ Chapoton, and R. Schiffler,
{\it Quivers with relations arising from clusters ($A_n$ case)}, Trans.\
Amer.\ Math.\ Soc.\ {\bf 358} (2006), 1347--1364. 

% \bibitem{EHS} K.\ Erdmann, T.\ Holm, N.\ Snashall,
% {\it Twisted bimodules and Hochschild cohomology for self-injective
% algebras of class $A_n$, II}, Algebr.\ Represent.\ Theory {\bf 5}
% (2002), 457--482.

% \bibitem{ES} K.\ Erdmann, A.\ Skowro\'nski, {\it The stable
% Calabi-Yau dimension of tame symmetric algebras}, J.\ Math.\ Soc.\ 
% Japan {\bf 58} (2006), 97--123.

% \bibitem{FR} S.\ Fomin, N.\ Reading, {\it 
% Generalized cluster complexes and Coxeter combinatorics},
% Int.\ Math.\ Res.\ Not.\ {\bf 44} (2005), 2709--2757.

% \bibitem{FZ} S.\ Fomin, A.\ Zelevinsky, 
% {\it Cluster algebras I. Foundations}, J.\ Amer.\ Math.\ Soc.
% {\bf 15} (2002), no.\ 2, 497--529.

\bibitem{FZ2}  S.\ Fomin and A.\ Zelevinsky, {\it Cluster algebras II:
    Finite type classification}, Invent.\ Math.\ {\bf 154} (2003),
  63--121. 

%\bibitem{Gabriel}  P.\ Gabriel, {\it Auslander-Reiten sequences and
%representation-finite algebras}, pp.\ 1--71 in ``Representation Theory
%I'' (Proceedings of the Workshop on the Present Trends in
%Representation Theory, Ottawa, 1979), edited by V.\ Dlab and P.\
%Gabriel, Lecture Notes in Math., Vol.\ 831, Springer, Berlin, 1980. 

%\bibitem{Happel2}  D.\ Happel, {\it Auslander-Reiten triangles in the
%derived categories of finite-dimensional algebras}, Proc.\ Amer.\
%Math.\ Soc.\ {\bf 112} (1991), 641--648.

% \bibitem{Happel}  D.\ Happel, {\it On the derived category of a finite
% dimensional algebra}, Comment.\ Math.\ Helv.\ {\bf 62} (1987),
% 339--389. 

%\bibitem{Happel3}  \bysame, ``Triangulated categories in the
%representation theory of finite dimensional algebras'', London Math.\
%Soc.\ Lecture Note Ser., Vol.\ 119, Cambridge University Press,
%Cambridge, 1988.

% \bibitem{HolmJorgensenDE} T.\ Holm, P.\ J{\o}rgensen, {\it Cluster
% categories, selfinjective algebras, and stable Calabi-Yau
% dimensions: types D and E}, preprint (2008), extended version of
% {\tt math.RT/0612451}.

%\bibitem{HJ-CYdim} T.\ Holm, P.\ J{\o}rgensen, {\it Stable Calabi-Yau
%dimension for finite type selfinjective algebras},
%preprint (2007). {\tt math.RT/0703488}. 

\bibitem{Iyama2}  O.\ Iyama, {\it Auslander correspondence}, Adv.\
  Math.\ {\bf 210} (2007), 51--82.

\bibitem{Iyama1}  O.\ Iyama, {\it Higher dimensional Auslander-Reiten
    theory on maximal orthogonal subcategories}, Adv.\ Math.\ {\bf
    210} (2007), 22--50.

\bibitem{IyamaYoshino}  O.\ Iyama and Y.\ Yoshino, {\it Mutation in
triangulated categories and rigid Cohen-Macaulay mo\-du\-les},
Invent.\ Math.\ {\bf 172} (2008), 117--168.

\bibitem{J}  P.\ J\o rgensen, {\it Auslander-Reiten theory over
topological spaces}, Comment.\ Math.\ Helv.\ {\bf 79} (2004),
160--182. 

%\bibitem{KellerDG}  B.\ Keller, {\it Deriving DG categories}, Ann.\
%Sci.\ \'Ecole Norm.\ Sup.\ (4) {\bf 27} (1994), 63--102.

\bibitem{Keller}  B.\ Keller, {\it On triangulated orbit categories},
Documenta Math.\ {\bf 10} (2005), 551--581.

%\bibitem{KellerReiten}  \bysame, {\it Acyclic Calabi-Yau categories
%are cluster categories}, in Oberwolfach Rep.\ {\bf 3}, no.\ 2 (2006).

\bibitem{KellerReiten2}  B.\ Keller and I.\ Reiten, {\it Cluster
tilted algebras are Gorenstein and stably Calabi-Yau},
Adv.\ Math.\ {\bf 211} (2007), 123-151.

\bibitem{KYZ}  B.\ Keller, D.\ Yang, and G.\ Zhou, The Hall algebra of
  a spherical object, preprint (2008).  {\tt math.RT/0810.5546}.

\bibitem{KoenigZhu} S.\ K\"onig and B.\ Zhu, {\it From triangulated
categories to abelian categories --- cluster tilting in a general
framework}, Math.\ Z.\ {\bf 258} (2008), 143--160.

\bibitem{K}  K.\ Kuribayashi, On the levels of spaces over a formal
space, preprint (2008).

%\bibitem{MacLane}  S.\ Mac Lane, ``Categories for the working
%mathematician'', Grad.\ Texts in Math., Vol.\ 5, Springer, Berlin,
%1971. 

% \bibitem{MiyachiYekutieli}  J.-I.\ Miyachi, A.\ Yekutieli, {\it
% Derived Picard groups of finite-dimensional hereditary algebras},
% Compositio Math.\ {\bf 129} (2001), 341--368.

%\bibitem{Neemanbook}  A.\ Neeman, ``Triangulated categories'', Ann.\
%of Math.\ Stud., Vol.\ 148, Princeton University Press, Princeton,
%2001. 

\bibitem{RVdB}  I.\ Reiten and M.\ Van den Bergh, {\it Noetherian
hereditary abelian categories sa\-tis\-fy\-ing Serre duality}, J.\
Amer.\ Math.\ Soc.\ {\bf 15} (2002), 295--366.

% \bibitem{Riedtmann}  C.\ Riedtmann, {\it Representation-finite
% self-injective algebras of class $A_{n}$}, pp.\ 449--520 in
% ``Representation theory, II'' (Proceedings of ICRA II, Ottawa, 1979),
% Lecture Notes in Math., Vol.\ 832, Springer, Berlin, 1980.

%\bibitem{RingTame}  C.\ M.\ Ringel, ``Tame algebras and quadratic
%forms'', Lecture Notes in Math., Vol. 1099, Springer, Berlin, 1984.

%\bibitem{Schiffler} R.\ Schiffler, 
%A geometric model for cluster categories of type $D_n$,
%preprint (2006). {\tt math.RT/0608264}.

% \bibitem{Skowronski} A.\ Skowro\'{n}ski, {\it Selfinjective algebras:
%     finite and tame type}, pp.\ 169--238 in: J.\ A.\ de la Pe\~{n}a,
%   R.\ Bautista, ``Trends in representation theory of algebras and
%   related topics'', Contemp.\ Math., Vol.\ 406, American Mathematical
%   Society, Rhode Island, 2006.
  
\bibitem{Thomas} H.\ Thomas, {\it Defining an $m$-cluster category},
J.\ Algebra {\bf 318} (2007), 37--46.
  
\bibitem{BinZhu} Bin Zhu, {\it Generalized cluster complexes via
quiver representations}, J.\ Algebraic Combin.\ {\bf 27} (2008),
25--54.

\end{thebibliography}
\end{document}